\documentclass[]{scrartcl}

\usepackage[utf8]{luainputenc}
\usepackage[USenglish]{babel}
\usepackage{csquotes}

\usepackage[a4paper,top=27mm,bottom=20mm,inner=25mm,outer=20mm]{geometry}

\usepackage[%
  backend=bibtex,bibencoding=ascii,
  style=numeric-comp,
  giveninits=true, uniquename=init, 
  natbib=true,
  url=true,
  doi=true,
  isbn=false,
  backref=false,
  maxnames=99,
  ]{biblatex}
\usepackage{amsmath}
\allowdisplaybreaks
\numberwithin{equation}{section}
\usepackage{amssymb}
\usepackage{commath}
\usepackage{mathtools}
\usepackage{bbm}
\usepackage{nicefrac}
\usepackage{subdepth}

\usepackage{siunitx}
\usepackage{algorithm}
\usepackage{algpseudocode}

\usepackage{amsthm}

\usepackage[plainpages=false,pdfpagelabels,hidelinks,unicode]{hyperref}

\usepackage{cleveref}

\usepackage{thmtools}
\usepackage{etoolbox}
\makeatletter
\patchcmd{\thmt@setheadstyle}
 {\bgroup\thmt@space}
 {\thmt@space}
 {}{}
\patchcmd{\thmt@setheadstyle}
 {\egroup\fi}
 {\fi}
 {}{}
\makeatother
\declaretheoremstyle[
  bodyfont=\normalfont\itshape,
  headformat=\NAME\ \NUMBER\NOTE,
]{myplain}
\declaretheoremstyle[
  headformat=\NAME\ \NUMBER\NOTE,
]{mydefinition}
\newcommand{\envqed}{{\lower-0.3ex\hbox{$\triangleleft$}}}
\declaretheorem[style=myplain,numberwithin=section]{theorem}

\declaretheorem[style=mydefinition,numberlike=theorem,qed=\envqed]{remark}

\usepackage{color}
\usepackage{graphicx}
\usepackage[small]{caption}
\usepackage{subcaption}

\begingroup\expandafter\expandafter\expandafter\endgroup
\expandafter\ifx\csname pdfsuppresswarningpagegroup\endcsname\relax
\else
\fi

\usepackage{booktabs}
\usepackage{rotating}
\usepackage{multirow}

\usepackage{enumitem}

\usepackage{ifluatex}
\ifluatex
  \usepackage[no-math]{fontspec}
\else
  \usepackage[T1]{fontenc}
\fi
\usepackage{newpxtext,newpxmath}

\let\epsilon\varepsilon
\let\phi\varphi
\let\rho\varrho

\newcommand{\N}{\mathbb{N}}

\newcommand{\I}{\operatorname{I}}

\renewcommand{\vec}[1]{\pmb{#1}}

\newcommand{\orcid}[1]{ORCID:~\href{https://orcid.org/#1}{#1}}
\usepackage{authblk}

\newenvironment{keywords}{\par\textbf{Key words.}}{\par}
\newenvironment{AMS}{\par\textbf{AMS subject classification.}}{\par}

\title{Convergence of hyperbolic approximations to higher-order PDEs for smooth solutions}

\author[1]{Jan Giesselmann\thanks{\orcid{0009-0008-0217-7244}}}
\affil[1]{Department of Mathematics, Technical University of Darmstadt, Germany}

\author[2]{Hendrik Ranocha\thanks{\orcid{0000-0002-3456-2277}}}
\affil[2]{Institute of Mathematics, Johannes Gutenberg University Mainz, Germany}

\date{January 2, 2026} 

\makeatletter
\hypersetup{pdfauthor={Jan Giesselmann, Hendrik Ranocha}} 
\hypersetup{pdftitle={Convergence of hyperbolic approximations to higher-order PDEs for smooth solutions}} 
\makeatother

\begin{document}

\maketitle

\begin{abstract}
\noindent
  We prove the convergence of hyperbolic approximations for several classes of higher-order PDEs, including the Benjamin-Bona-Mahony, Korteweg-de Vries, Gardner, Kawahara, and Kuramoto-Sivashinsky equations, provided a smooth solution of the limiting problem exists.
We only require weak (entropy) solutions of the hyperbolic approximations.
Thereby, we provide a solid foundation for these approximations, which have been used in the literature without rigorous convergence analysis.
We also present numerical results that support our theoretical findings.

\end{abstract}



\section{Introduction}

Hyperbolic approximations, also known as hyperbolic relaxations or hyperbolizations, of higher-order partial differential equations (PDEs) have been studied for a long time, such as the hyperbolic heat equation discussed by Cattaneo and Vernotte \cite{cattaneo1958forme,vernotte1958paradoxes}.
These hyperbolic approximations have been applied in various contexts, including elliptic problems
\cite{nishikawa2007first,ruter2018hyperbolic,schlottkelakemper2021purely,ranocha2023discontinuous},
dispersive wave equations
\cite{gavrilyuk2022hyperbolic,besse2022perfectly,busto2021high,favrie2017rapid,ranocha2025structure},
and other higher-order PDEs
\cite{dhaouadi2024first,gavrilyuk2024conduit,ketcheson2025approximation}.
Typically, such hyperbolizations are developed for specific PDEs using tailored techniques with the goal of maintaining essential structures of the original PDE.
Notably, an augmented Lagrangian approach \cite{favrie2017rapid,Dhaouadi2019,Favire2024} has proven successful.
Ketcheson and Biswas \cite{ketcheson2025approximation} introduced a broad framework to create stable hyperbolizations of a class of linear PDEs and presented numerical results suggesting that this method may also be effective for nonlinear and more complex PDEs.
They additionally offer a review of several hyperbolizations along with motivations for developing them.

As far as we know, hyperbolizations have been used quite extensively, but usually without a rigorous convergence analysis to support their approximation properties.
Some notable exceptions include the hyperbolic heat equation \cite{nagy1994behavior},
and the Serre-Green-Naghdi system at the analytical level	\cite{duchene2019rigorous},
as well as the Korteweg-de Vries and Benjamin-Bona-Mahony equations at the numerical level \cite{biswas2025traveling,bleecke2025asymptotic}.

In this paper, we provide a rigorous convergence analysis of a broad class of hyperbolic approximations to higher-order PDEs using the relative energy/entropy method.
We start with a structure similar to the Benjamin-Bona-Mahony equation that involves mixed space-time derivatives in Section~\ref{Sec:mst}.
We introduce the relative energy method in this context and explain some choices we make in the convergence analysis.
Section~\ref{sec:spatial} is dedicated to more general classes of PDEs with higher-order spatial derivatives.
We perform the convergence analysis for each class and present numerical results that support our findings.
Finally, we summarize our results in Section~\ref{sec:summary} and discuss potential extensions of our work.

The relative entropy (sometimes also called relative energy) method has a long history for both hyperbolic and parabolic problems.
It relies on leveraging an energy or entropy structure of the considered system in order to control the distance between two solutions.
Depending on the context, the functionals of interest are either called energy or entropy.
In this paper, we will use the term relative energy consistently (although we could call it relative entropy equally well). For consistency with the literature, we will still use the common notion of weak entropy solutions (and will not call them weak energy solutions).

For the use of the relative energy method in parabolic models, in particular with respect to studying the convergence to steady states, we refer to \cite{Juengel2016} and references therein.
For hyperbolic models, it dates back to Dafermos \cite{Dafermos1979} and di~Perna \cite{Diperna1979}, where it was, in particular, used to establish weak-strong uniqueness.
In addition to being a well-established tool for investigating the stability of solutions to systems of hyperbolic conservation laws, it also has a long history as a tool for describing the relationship between solutions of different models, where usually one can be understood as the relaxation limit of the other.
It has been particularly successful in identifying low Mach and large friction limits of hyperbolic systems, see \cite{Tzavaras2005,Lattanzio2017,GLT2017,Feireisl2024,Egger2023,Gallenmueller2024}, but also in comparing different moment approximations of kinetic equations \cite{Alldredge2023}.
For an extension to hyperbolic-parabolic systems, see \cite{Christoforou2018}.
A general property in this analysis is that the solution to the limiting system needs to be a strong solution whereas the solutions to the approximating system can be allowed to be (weak) entropy solutions.
In this paper, we leverage the fact that many hyperbolic approximations of higher-order equations, as derived in \cite{ketcheson2025approximation}, are equipped with a (rather simple) energy that allows for relative energy estimates.
One feature that creates challenges not encountered in previous works is that the energy degenerates, i.e., the modulus of convexity approaches zero in many directions when the limit problem is approached.
We explain how we handle this in Section~\ref{ssec:bica}.
Let us also notice that several hyperbolic approximations suggested in the literature \cite{ketcheson2025approximation} do not admit a (simple) energy.
In these cases, we suggest alternative hyperbolic approximations that do admit a (simple) energy.
Incidentally, these new approximations are easier to handle numerically in the sense that we can construct energy-preserving discretizations.

To demonstrate the analytical convergence results, we implement several model problems
numerically using (upwind) summation-by-parts (SBP) operators in space and
(additive) Runge-Kutta methods in time. More details on upwind SBP operators can be
found in \cite{mattsson2017diagonal,ranocha2021broad}; applications to hyperbolic
approximations of the KdV and BBM equations are discussed in
\cite{biswas2025traveling,bleecke2025asymptotic}. A more general discussion of the
background of SBP operators is contained in the review articles
\cite{svard2014review,fernandez2014review}. Although we mainly focus on finite
differences \cite{kreiss1974finite,strand1994summation,carpenter1994time},
the general framework of SBP operators also includes
finite volumes \cite{nordstrom2001finite},
finite elements \cite{hicken2016multidimensional,hicken2020entropy,abgrall2020analysisI},
discontinuous Galerkin~(DG) methods \cite{gassner2013skew,carpenter2014entropy},
flux reconstruction~(FR) \cite{huynh2007flux,vincent2011newclass,ranocha2016summation},
active flux methods \cite{eymann2011active,barsukow2025stability},
as well as meshless schemes \cite{hicken2024constructing}.
We use the time integration schemes of \cite{ascher1997implicit} and refer to
\cite{boscarino2024asymptotic,boscarino2024implicit} for background information on
implicit-explicit (IMEX) and asymptotic-preserving methods.

All source code and data required to reproduce the numerical results are available
in our reproducibility repository \cite{giesselmann2025convergenceRepro}. We implemented
all methods in Julia \cite{bezanson2017julia}, used SummationByPartsOperators.jl
\cite{ranocha2021sbp} for the spatial discretizations, and created the figures
using Makie.jl \cite{danisch2021makie}.

\section{Dispersive equations with mixed space-time derivatives}\label{Sec:mst}

First, we concentrate on PDEs of the form
\begin{equation}
\label{eq:mixed_limit}
  \partial_t u(t, x) + \partial_x f(u(t,x)) - \partial_x^2 \partial_t u(t, x) = 0
\end{equation}
with periodic boundary conditions and appropriate initial data.
A classical example is the Benjamin-Bona-Mahony (BBM) equation \cite{benjamin1972model}
\begin{equation}
\label{eq:bbm}
  \partial_t u(t, x) + \partial_x \frac{u(t,x)^2}{2} - \partial_x^2 \partial_t u(t, x) = 0,
\end{equation}
where $f(u) = u^2 / 2$. Note that the typical form of the BBM equation
includes an additional linear term $+ \partial_x u(t, x)$ that is not
present in \eqref{eq:bbm} since we follow the notation of \cite{gavrilyuk2022hyperbolic}.
The typical form can be obtained by a change of variables $u \mapsto u + 1$.

Gavrilyuk and Shyue \cite{gavrilyuk2022hyperbolic} introduced a hyperbolic
approximation of the BBM equation \eqref{eq:bbm} based on an augmented Lagrangian.
Generalizing their system to a general flux $f(u)$ in \eqref{eq:mixed_limit}
yields\footnote{We use the slightly different ordering $q_0 = u$, $q_1 = w$, and $q_2 = v$ of the equations compared to the notation of \cite{gavrilyuk2022hyperbolic} for better consistency with the analysis presented in Section~\ref{sec:spatial}. We also use only a single hyperbolic relaxation parameter $c = \widehat{c} = 1 / \tau$.}
\begin{equation}
\label{eq:mixed_hyperbolic}
\begin{aligned}
  \partial_t q_{0}
  + \partial_x f(q_{0})
  + \partial_x q_{2}
  &= 0,
  \\
  \partial_t q_{1}
  + \tau \partial_x q_{1}
  &= -q_{2},
  \\
  \tau \partial_t q_{2}
  + \partial_x q_{0}
  &= q_{1},
\end{aligned}
\end{equation}
where $\tau > 0$ is the hyperbolic relaxation parameter. In the limit $\tau \to 0$,
we formally obtain $q_1 \to \partial_x q_0$ from the third equation and
$q_2 \to -\partial_{t} q_1 = -\partial_{tx} q_0$ from the second equation,
so that $q_0 \to u$ formally satisfies the limit equation \eqref{eq:mixed_limit}.

Similar to \cite{gavrilyuk2022hyperbolic}, the system \eqref{eq:mixed_hyperbolic}
is hyperbolic with Jacobian eigenvalues
\begin{equation*}
  \tau \quad\text{and}\quad \frac{1}{2} f'(q_0) \pm \sqrt{(f'(q_0))^2 + 4 \tau^{-1}}
\end{equation*}
and conserves the energy
\begin{equation}
  \eta(q) = \int \left( \frac{1}{2} q_0^2 + \frac{1}{2} q_1^2 + \frac{\tau}{2} q_2^2 \right) \dif x.
\end{equation}

\subsection{Basic idea of the convergence analysis}\label{ssec:bica}

The basic idea of the convergence analysis is to interpret a sufficiently regular solution of the limiting equation as the solution of a perturbed version of the approximating system.
Here, perturbed means that there are additional \emph{residual} terms that scale with $\tau$.
Then, one can use the relative energy stability framework for the approximate system to bound the difference between the exact solution of the approximating system and the solution to the limiting system in terms of these residuals.
Since the residuals scale with $\tau$ and otherwise only depend on the limiting solution and its derivatives, this allows us to infer convergence.
This is very similar to the procedure used in, e.g., \cite{Lattanzio2017,Egger2023}. A particular challenge here is that the energy of the approximate system degenerates for $\tau \rightarrow 0$ in the sense that in many (but fortunately not all) directions the modulus of convexity goes to zero for $\tau \rightarrow 0$.
This leads to a situation where residuals in the first and second equations are easier to handle than a residual in the third equation.
Thus, if $u$ denotes the solution to the limiting problem we do not use the seemingly obvious definition $\bar q=(u, \partial_x u, -\partial_{tx} u)$ to create an approximate solution to the approximating system, but we perturb this by terms of order $\tau$ in such a way that there is only a residual in the second equation, while the other equations are satisfied exactly.
If we used $\bar q=(u, \partial_x u, -\partial_{tx} u)$, we would have a residual in the third equation and our analysis would only provide convergence with rate $\sqrt{\tau}$.

\subsection{Convergence analysis}

Suppose $f \in W^{2,\infty}(\mathbb{R})$ and we are given a solution $u$ to \eqref{eq:mixed_limit} and suppose $u \in H^4((0,T) \times \Omega)$ and
$\partial_x u \in L^\infty((0,T)\times \Omega)$,
where $\Omega$ is our computational domain (we always assume periodic boundary conditions).
Then, we define an approximate solution
\begin{equation}
  \bar{q} = (\bar q_0, \bar q_1, \bar q_2) := (u, \partial_x u- \tau \partial_{ttx} u, -\partial_{tx} u )
\end{equation}
to \eqref{eq:mixed_hyperbolic}.
Inserting $\bar q$ into the hyperbolic approximation \eqref{eq:mixed_hyperbolic}, we obtain
\begin{equation}
\label{eq:approximation_mixed_res}
\begin{aligned}
  \partial_t \bar q_{0}
  + \partial_x f(\bar q_{0})
  + \partial_x \bar q_{2}
  &= 0,
  \\
  \partial_t \bar q_{1}
  + \tau \partial_x \bar q_{1}
  &= -\bar q_{2} + \tau R, 
  \\
  \tau \partial_t \bar q_{2}
  + \partial_x \bar q_{0}
  &= \bar  q_{1},
\end{aligned}
\end{equation}
where the residual is given by $R = -\partial_{tttx}u + \partial_{xx} u - \tau \partial_{ttxx} u$.

Let us now consider the time evolution of the relative energy
\begin{equation}
 \eta(q, \bar q) = \int_\Omega \left(
    \frac12 (q_0 - \bar q_0)^2 + \frac12 (q_1 - \bar q_1)^2 + \frac{\tau}{2} (q_2 - \bar q_2)^2
   \right) \dif x.
\end{equation}
We obtain
\begin{equation}
\begin{aligned}
  \frac{\dif}{\dif t} \eta(q, \bar q)
  &=
  \int_\Omega \left(
    (q_0 - \bar q_0) \partial_t (q_0- \bar q_0)
    + (q_1 - \bar q_1) \partial_t (q_1- \bar q_1)
    + \tau (q_2 - \bar q_2) \partial_t (q_2- \bar q_2)
  \right) \dif x
  \\
  &=
  \int_\Omega \biggl(
    - (q_0 - \bar q_0) \partial_x \bigl(f(q_0)- f(\bar q_0)\bigr)
    - (q_0 - \bar q_0) \partial_x \bigl(q_2 - \bar q_2\bigr)
    \\
    &\qquad\quad
    - (q_1 - \bar q_1)\tau \partial_x (q_1- \bar q_1) - (q_1 - \bar q_1) (q_2 - \bar q_2)
    \\
    &\qquad\quad
    - (q_1 - \bar q_1) \tau R -  (q_2 - \bar q_2)\partial_x (q_0- \bar q_0) +  (q_2 - \bar q_2) (q_1 - \bar q_1)
  \biggr) \dif x
  \\
  &=
  -\int_\Omega \left(
    (q_0 - \bar q_0)\partial_x (f(q_0)- f(\bar q_0))  + (q_1 - \bar q_1)\tau R
  \right) \dif x.
\end{aligned}
\end{equation}
Let us define the entropy flux $\mu$ such that $\mu'(a) = a f'(a)$ for all $a \in \mathbb{R}$.
Then,
\begin{multline}
 (q_0 - \bar q_0)\partial_x (f(q_0)- f(\bar q_0)) \\
 =
 \partial_x \bigl(\mu(q_0) - \mu(\bar q_0) - \bar q_0 (f(q_0) - f(\bar q_0))\bigr)
 + (\partial_x \bar q_0) \bigl(f(q_0) - f(\bar q_0) - f'(\bar q_0)(q_0 - \bar q_0)\bigr).
\end{multline}
Inserting this into the evolution equation for $\eta(q, \bar q)$ leads to
\begin{equation}
  \frac{\dif}{\dif t} \eta(q, \bar q)
  =
  - \int_\Omega \left(
    (\partial_x \bar q_0) \bigl(f(q_0) - f(\bar q_0) - f'(\bar q_0)(q_0 - \bar q_0)\bigr) + (q_1 - \bar q_1)\tau R
  \right) \dif x.
\end{equation}
Then, Young's inequality implies
\begin{equation}
  \frac{\dif}{\dif t} \eta(q, \bar q)
  \leq
  (\|\partial_x u\|_{L^\infty} +1 ) \eta(q, \bar q)
  + \tau^2 \|\partial_{tttx} u\|_{L^2}^2 + \tau^2\|\partial_{xx} u\|_{L^2}^2 + \tau^4 \|\partial_{ttxx} u\|_{L^2}^2.
\end{equation}
We can conclude
\begin{equation}
  \| \eta(q, \bar q)\|_{L^\infty(0,T)} = \mathcal{O}(\tau^2)
\end{equation}
by using Gronwall's lemma since $u$ is independent of $\tau$
and $\eta(q, \bar q)|_{t=0}= \mathcal{O}(\tau^2)$ if we initialize the hyperbolic
approximation \eqref{eq:mixed_hyperbolic} with
$q_0 = u$, $q_1 = \partial_x u$, and $q_2 = -\partial_{tx} u$.
Summing up, we have proved the following
\begin{theorem}
\label{thm:convergence_mixed}
  Let $T>0$ and $f \in W^{2,\infty}_\mathrm{loc}(\mathbb{R})$ such that $f'' \in L^\infty(\mathbb{R})$.
  Let $u \in H^{4}((0,T) \times \Omega)$ such that $\partial_x u \in L^\infty((0,T)\times \Omega)$ be a solution to \eqref{eq:mixed_limit} with initial data $u_0 \in H^{4}(\Omega)$.
  Let, for each $\tau>0$, $q$ be an entropy solution to \eqref{eq:mixed_hyperbolic} with $q|_{t=0}= (u_0, \partial_x u_0, -\partial_{tx} u|_{t=0} )$. Then
  \begin{equation}
    \| u - q_0 \|_{L^\infty(0,T, L^2(\Omega))}
    + \| \partial_x u - q_1 \|_{L^\infty(0,T, L^2(\Omega))}
    = \mathcal{O}(\tau).
  \end{equation}
\end{theorem}

\begin{remark}
\label{rem:convergence_mixed}
  It should be noted that $f'' \in L^{\infty}(\mathbb{R})$ is not strictly needed:
  What we need is that there exists $C>0$ such that
  \begin{equation*}
    |f(q_0) - f(\bar q_0) - f'(\bar q_0)(q_0 - \bar q_0)|
    \leq
    C |q_0 - \bar q_0|^2 \quad \text{ uniformly in } (0,T) \times \Omega.
  \end{equation*}
  This is also satisfied if $f \in C^2(\mathbb{R})$ and if there exists some compact $K \subset \mathbb{R}$ such that $\bar q_0$ (as constructed from $u$) and $q_0$ (for all values of $\tau$) only take values in this set.
\end{remark}

\subsection{Numerical demonstration}

We use the numerical methods of \cite{bleecke2025asymptotic} to check the
convergence of the hyperbolic approximation \eqref{eq:mixed_hyperbolic}
to the BBM equation \eqref{eq:bbm} numerically. Thus, we discretize the BBM
equation \eqref{eq:bbm} as
\begin{equation}
  \partial_t \vec{u}
  + \frac{1}{3} (\I - D_+ D_-)^{-1} \left(
    \vec{u} D_1 \vec{u} + D_1 \vec{u}^2
  \right)
  =
  \vec{0},
\end{equation}
where $D_\pm$ are periodic upwind SBP (summation by parts) operators and
$D_1 = (D_+ + D_-) / 2$ the corresponding central SBP operator. We apply the
explicit part of ARS(4,4,3) \cite{ascher1997implicit} to the resulting system.
Similarly, we discretize the hyperbolic approximation \eqref{eq:mixed_hyperbolic} as
\begin{equation}
\label{eq:bbmh_semidiscretization}
  \partial_t
  \begin{pmatrix}
    \vec{q_0} \\
    \vec{q_1} \\
    \vec{q_2}
  \end{pmatrix}
  + \begin{pmatrix}
      \frac{1}{3} (\vec{q_0} D_1 \vec{q_0} + D_1 \vec{q_0}^{\!\! 2}) \\
      \tau D_1 \vec{q_1} \\
      \vec{0}
  \end{pmatrix}
  + \begin{pmatrix}
    D_1 \vec{q_2} \\
    \vec{q_2} \\
    \tau^{-1} D_1 \vec{q_0} - \tau^{-1} \vec{q_1}
  \end{pmatrix}
  =
  \begin{pmatrix}
    \vec{0} \\
    \vec{0} \\
    \vec{0}
  \end{pmatrix},
\end{equation}
where we apply the explicit part of the time integration scheme to the first
term and the implicit part to the second term. Please note that we used a split
form of the nonlinear term to guarantee energy conservation since the chain rule
is not available in general discretely \cite{ranocha2019mimetic}.
The semidiscretization conserves discrete analogs of $\int q_0 \dif x$ and the energy
$\int (q_0^2 + q_1^2 + \tau q_2^2) \dif x / 2$.
An analysis of the structure- and asymptotic-preserving properties of these
discretizations is contained in \cite{bleecke2025asymptotic}.

\begin{figure}[htb]
  \includegraphics[width=\textwidth]{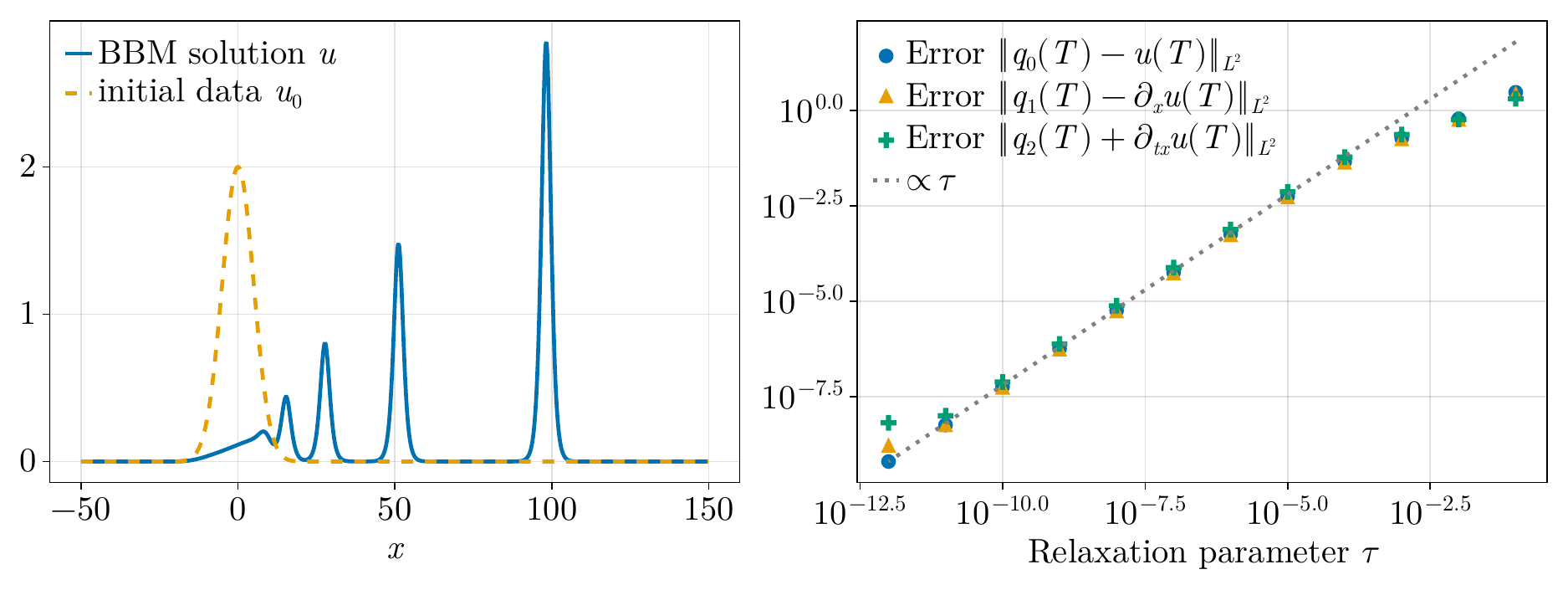}
  \caption{Convergence of the hyperbolic approximation \eqref{eq:bbmh_semidiscretization}
           to the BBM equation \eqref{eq:bbm}. The left plot shows the numerical solution
           at the final time $T = 100$ and the initial condition. The right plot shows
           the convergence of the discrete $L^2$ error at the final time $T$
           as a function of $\tau$.}
  \label{fig:benjamin_bona_mahony_convergence}
\end{figure}

To test convergence of the hyperbolic approximation, we use the initial data
\begin{equation}
  u(0, x) = 2 \exp\bigl(-0.02 x^2\bigr)
\end{equation}
in the domain $\Omega = [-50, 150]$ with periodic boundary conditions. The numerical
solution of the BBM equation as well as the discrete $L^2$ errors are shown in
Figure~\ref{fig:benjamin_bona_mahony_convergence}
for $2^{10}$ grid points, seventh-order accurate upwind operators, and
time step size $\Delta t = 0.1$.
As expected, the Gaussian initial condition splits into traveling waves at
the final time $t = 100$. Moreover, the hyperbolic approximation converges
to the BBM solution with the expected order $O(\tau)$ predicted by
Theorem~\ref{thm:convergence_mixed}.
Notably, the derivative approximations converge with the same order instead of
a reduced order as one might expect from the analysis in Theorem~\ref{thm:convergence_mixed}.

\section{Dispersive equations with purely spatial higher-order derivatives}
\label{sec:spatial}

In this section, we consider PDEs with purely spatial higher-order derivatives.
Ketcheson and Biswas \cite{ketcheson2025approximation} analyzed linear equations
of the form
\begin{equation}
\label{eq:ketcheson_biswas_limit}
  \partial_t u(t, x) + \sigma_0 \partial_x^m u(t, x) = 0
\end{equation}
with $\sigma_0 \in \{-1, +1\}$, $m \in \N$, and periodic boundary conditions.
They constructed stable hyperbolic approximations\footnote{Please note the
misprints in \cite{ketcheson2025approximation}. The statements of the
permutation matrix $P$ is correct, but the resulting equation (27)
contains some errors in the original journal publication
\cite{ketcheson2025approximation}. We thank David Ketcheson for helping
us to clarify this issue.}
\begin{equation}
\begin{aligned}
\label{eq:ketcheson_biswas_hyperbolic}
  \partial_t q_{0}
  + \sigma_0 \partial_x q_{m-1}
  &= 0,
  \\
  \tau \partial_t q_{j}
  + \sigma_0 (-1)^{j} \partial_x q_{m-j-1}
  &= \sigma_0 (-1)^{j} q_{m-j},
  & j &\in \{1, \dots, \lceil m / 2 \rceil - 1\},
  \\
  \tau \partial_t q_{j}
  + \sigma_0 (-1)^{j+m-1} \partial_x q_{m-j-1}
  &= \sigma_0 (-1)^{j+m-1} q_{m-j},
  & j &\in \{\lceil m / 2 \rceil, \dots, m-1\},
\end{aligned}
\end{equation}
where $\tau > 0$ is the parameter of the hyperbolization,
$q_0$ is the main variable that formally converges to $u$ in the limit
$\tau \to 0$, and the other $q_j$ are auxiliary variables that formally
converge to the $j$-th spatial derivative $\partial_x^j u$.

Ketcheson and Biswas \cite{ketcheson2025approximation} showed that
the system \eqref{eq:ketcheson_biswas_hyperbolic} has the quadratic energy
\begin{equation}
\label{eq:spatial_energy}
  \eta(q) = \int \left( \frac{1}{2} q_0^2 + \sum_{j=1}^{m-1} \frac{\tau}{2} q_j^2 \right) \dif x.
\end{equation}
Concretely, the system \eqref{eq:ketcheson_biswas_hyperbolic} conserves
$\eta(q)$  for odd $m$
and dissipates $\eta(q)$ for even $m$ with $\sigma_0 = (-1)^{m/2}$.
Indeed, for odd $m$, smooth solutions of \eqref{eq:ketcheson_biswas_hyperbolic} satisfy
\begin{equation}
\frac{\dif}{\dif t} \eta(q)=0
\end{equation}
and for even $m$, smooth solutions of \eqref{eq:ketcheson_biswas_hyperbolic} satisfy
 \begin{equation}
   \frac{\dif}{\dif t}\eta(q)
   = -\sigma_0 (-1)^{m/2} \int q_{m/2}^2 \dif x
   \le 0.
 \end{equation}
Both were already shown in \cite{ketcheson2025approximation} and also follow from the relative energy computations later in this paper by setting $\bar q_0 = \dots = \bar q_m = 0$ and $R=0$.

Next, we generalize these hyperbolic approximations by allowing for an additional
nonlinear first-order term $+ \partial_x f(u)$ in \eqref{eq:ketcheson_biswas_limit}
as well as additional higher-order
spatial derivatives. Since the analysis depends strongly on the parity of the leading
order $m$, we distinguish between the cases of odd and even $m$ in the following.

\subsection{Odd leading order \texorpdfstring{$m$}{m}}

We first consider partial differential equations (PDEs) of the form
\begin{equation}
\label{eq:spatial_odd_limit}
  \partial_t u(t, x)
  + \partial_x f(u(t,x))
  - \mu \partial_{xx} u(t,x)
  + \sigma_0 \partial_x^m u(t, x) = 0
\end{equation}
with $\sigma_0 \in \{-1, +1\}$, $\mu \ge 0$, and odd $m \in \N$. We generalize
the hyperbolic approximations \eqref{eq:ketcheson_biswas_hyperbolic}
of \cite{ketcheson2025approximation} to
\begin{equation}
\label{eq:spatial_odd_hyperbolic}
\begin{aligned}
  \partial_t q_{0}
  + \partial_x f(q_{0})
  + \sigma_0 \partial_x q_{m-1}
  &= 0,
  \\
  \tau \partial_t q_{j}
  + \sigma_0 (-1)^{j} \partial_x q_{m-j-1}
  &= \sigma_0 (-1)^{j} q_{m-j} - \delta_{j,1} \mu q_1,
  & j &\in \{1, \dots, m - 1\},
\end{aligned}
\end{equation}
where $\delta_{j,1}$ is the Kronecker delta. Please note the specific choice
we made to approximate the dissipative term $-\mu \partial_{xx} u$ by introducing
an additional term in the evolution equation for $q_1$. This choice ensures that
the hyperbolic approximation conserves the quadratic energy \eqref{eq:spatial_energy}
for $\mu = 0$ and results in energy dissipation
\begin{equation}
  \frac{\dif}{\dif t} \eta(q)
  =
  -\mu \int q_1^2 \dif x
  \leq 0
\end{equation}
for $\mu \ge 0$.

\begin{remark}
 The philosophy of the convergence analysis here and in what follows is very similar to the one in Section~\ref{Sec:mst}. There is one difference:
 Here the structure of the energy is such that a residual in the first equation is easier to handle than a  residual in any of the other equations. With the seemingly obvious definition $\bar q=(u, \partial_x u, \dots,\partial_{x}^{m-2} u, \partial_{x}^{m-1} u - \frac{\mu}{\sigma_0} \partial_x u )$ for an approximate solution to the approximating system, we would have residuals in all but the first equation which would lead to  convergence with rate $\sqrt{\tau}$. Instead, we perturb this by terms of order $\tau$ in such a way that there is only a residual in the first equation, while the other equations are satisfied exactly.
\end{remark}

To perform the convergence analysis using the relative energy method,
we assume to be given a smooth solution $u$ to \eqref{eq:spatial_odd_limit}.
As explained above, we wish to define $\bar q$ such that it satisfies
\begin{equation}
\label{eq:spatial_odd_res}
\begin{aligned}
  \partial_t \bar q_{0}
  + \partial_x f(\bar q_{0})
  + \sigma_0 \partial_x \bar q_{m-1}
  &= \tau R,
  \\
  \tau \partial_t \bar q_{j}
  + \sigma_0 (-1)^{j} \partial_x \bar q_{m-j-1}
  &= \sigma_0 (-1)^{j} \bar q_{m-j} - \delta_{j,1} \sigma_0 \mu \bar q_1,
  & j &\in \{1, \dots, m - 1\},
\end{aligned}
\end{equation}
where the residual $R$ depends on $u$ and its derivatives and is uniformly
bounded for $ \tau \to 0$. The crucial point is that there is only a residual
in the first equation.

To define $\bar q$ we begin by setting $q_{\frac{m-1}{2}}= \partial_x^{\frac{m-1}{2}} u$,
and then we do the following two steps successively for all odd $k \in \{1, \dots, m-2\}$:
\begin{enumerate}
  \item We use the evolution equation for $q_j$ with $j=\frac{m-k}{2}$:
        \begin{equation}
          \tau \partial_t q_{\frac{m-k}{2}} + \sigma_0 (-1)^{\frac{m-k}{2}} \partial_x q_{\frac{m+k-2}{2}} = \sigma_0 (-1)^{\frac{m-k}{2}} q_{\frac{m+k}{2}}
        \end{equation}
        to determine $q_{\frac{m+k}{2}}$.
        Note that $q_{\frac{m-k}{2}}$ and $q_{\frac{m+k-2}{2}}$ have already been
        determined in the iteration.
  \item We use the evolution equation for $q_j$ with $j=\frac{m+k}{2}$:
        \begin{equation}
          \tau \partial_t q_{\frac{m+k}{2}} + \sigma_0 (-1)^\frac{m+k}{2} \partial_x q_{\frac{m-k-2}{2}} = \sigma_0 (-1)^\frac{m+k}{2} q_{\frac{m-k}{2}}
        \end{equation}
        to determine $q_{\frac{m-k-2}{2}}$.
        Note that $q_{\frac{m+k}{2}}$ and $q_{\frac{m-k}{2}}$ and each term
        in them contains at least one $x$-derivative.
        Strictly speaking $q_{\frac{m-k-2}{2}}$ is only defined up to a constant, but we can easily fix this constant such that the leading order term is $q_{\frac{m-k-2}{2}}=\partial_x ^{\frac{m-k-2}{2}} u.$
\end{enumerate}
By this procedure, $\bar q$ satisfies each equation in \eqref{eq:spatial_odd_hyperbolic},
except the first one, exactly.
In addition, for $j \in \{0, \dots, m-2\}$, each $\bar q_j$ satisfies
\begin{equation}\label{eq:approxbarq}
  \bar q_j= \partial_x^j u + \tau R_j,
\end{equation}
where $R_j$ depends on $u$ and its derivatives
and is uniformly bounded for $\tau \to 0$.
Similarly, $\bar q_{m-1}$ satisfies
\begin{equation}
  \sigma_0 \bar q_{m-1}= \sigma_0 \partial_x^{m-1} u - \mu \partial_x u + \tau R_{m-1},
\end{equation}
where $R_{m-1}$ depends on $u$ and its derivatives and is uniformly bounded for
$ \tau \to 0$.
To be more precise, if we write $j = \tfrac{m - 1 \pm k}{2}$ with $k$ even, then
$R_j$ only depends on derivatives of $u$ up to order $\frac{m-1}{2} + \frac{k}{2}$.
This is a consequence of the iterative definition of the $\bar q_j$:
initially $\bar q_{\frac{m-1}{2}}$ only depends on a derivative of order $\frac{m-1}{2}$.
In step one of the iteration, $\bar q_{\frac{m+k}{2}}$ will contain one
order of derivatives more than $\bar q_{\frac{m-k}{2}}$ and $\bar q_{\frac{m+k}{2}}$
contained, whereas in step two, $\bar q_{\frac{m-k-2}{2}}$ will contain derivatives
of at most the same order as those appearing in $\bar q_{\frac{m+k}{2}}$ and
$\bar q_{\frac{m-k}{2}}$.
The residual $R$ in \eqref{eq:spatial_odd_res} satisfies
\[R= \partial_t R_0 +  \frac{1}{\tau}\partial_x (f(\bar q_0) - f(u) ) + \sigma_0 \partial_x R_{m-1}.\]
We can infer that $R$ is uniformly bounded for small $\tau$ by using \eqref{eq:approxbarq}, Taylors formula and the product rule since the first two derivatives of $f$ are bounded on the (bounded) images of $\bar q_0$ and $u$.

We define the relative energy
\begin{equation}
  \eta(q, \bar q)
  :=
  \int_\Omega \left(
    \frac12 (q_0 - \bar q_0)^2 + \frac{\tau}{2} \sum_{i=1}^{m-1} (q_i - \bar q_i)^2
  \right) \dif x.
\end{equation}
Based on \eqref{eq:spatial_odd_limit} and \eqref{eq:spatial_odd_res}, we compute
\begin{equation}
\begin{aligned}
  \frac{\dif}{\dif t} \eta(q, \bar q)
  &=
  \int \left( (q_0 - \bar q_0) \partial_t (q_0 - \bar q_0) + \sum_{j=1}^{m-1} \tau (q_j - \bar q_j) \partial_t (q_j - \bar q_j) \right) \dif x
  \\
  &= -\int \biggl(
    (q_0 - \bar q_0) \partial_x \big(f(q_0) - f(\bar q_0)\big) + (q_0 - \bar q_0) \tau R+ \sigma_0 (q_0 - \bar q_0) \partial_x (q_{m-1}- \bar q_{m-1})
    \\
    &\qquad
    - \mu (q_1 - \bar q_1)^2+ \sum_{j=1}^{m - 1} \sigma_0 (-1)^j (q_{j} - \bar q_j)
      \bigl( \partial_x q_{m-j-1} - q_{m-j}- \partial_x \bar q_{m-j-1} +\bar  q_{m-j} \bigr)
  \biggr) \dif x
  \\
  &= -\int \biggl(
    (q_0 - \bar q_0) \partial_x \big(f(q_0) - f(\bar q_0)\big)
    + (q_0 - \bar q_0) \tau R
    - \mu (q_1 - \bar q_1)^2
    \\
    &\qquad
    + \sigma_0 \sum_{j=0}^{m - 1} (-1)^j (q_{j}- \bar q_j) \partial_x ( q_{m-j-1}- \bar q_{m-j-1})
    \\
    &\qquad
    - \sigma_0 \sum_{j=1}^{m - 1} (-1)^j (q_{j}- \bar q_j) (q_{m-j}- \bar q_{m-j})
  \biggr) \dif x
  \\
  &
  \leq \|\partial_x  \bar q_0 \|_{L^\infty}\eta(q,\bar q) + \tau^2 \int R^2 \dif x.
\end{aligned}
\end{equation}
We can conclude
\begin{equation}
  \| \eta(q, \bar q)\|_{L^\infty(0,T)} = \mathcal{O}(\tau^2)
\end{equation}
by using Gronwall's lemma since $u$ is independent of $\tau$.

Note that the definition of $\bar q_0$ contains derivatives of $u$
of degree $m-1$ so that the convergence proof requires
$u \in W^{m,\infty}((0,T) \times \Omega)$.
The preceding computation implies the following theorem
\begin{theorem}
\label{thm:convergence_spatial_odd}
  Let $T>0$ and $f \in W^{2,\infty}_{loc}(\mathbb{R})$ such that
  $f'' \in L^\infty(\mathbb{R})$.
  Let $u \in W^{m,\infty}((0,T) \times \Omega)$ be a solution to
  \eqref{eq:spatial_odd_limit} with initial data $u_0 \in W^{m,\infty}(\Omega)$.
  Let, for each $\tau>0$, $q$ be an entropy solution to \eqref{eq:spatial_odd_hyperbolic}
  with $q|_{t=0}= (u_0, \partial_x u_0 ,\dots, \partial_x^{m-1} u_0 )$. Then
  \begin{multline}
    \| u - q_0 \|_{L^\infty(0,T, L^2(\Omega))}
    + \sum_{j=1}^{m-2} \sqrt{\tau} \| \partial_x^j u - q_j \|_{L^\infty(0,T, L^2(\Omega))}
    \\
    + \sqrt{\tau} \|\sigma_0 \partial_x^{m-1} u -\mu \partial_x u  -\sigma_0 q_{m-1} \|_{L^\infty(0,T, L^2(\Omega))}
    =
    \mathcal{O}(\tau).
  \end{multline}
\end{theorem}

The result also holds in a more general setting similar to
Remark~\ref{rem:convergence_mixed}.
Next, we demonstrate the convergence results numerically for several examples.

\subsubsection{Korteweg-de Vries equation}

We apply the structure-preserving methods of \cite{biswas2025traveling} to check
the convergence of the hyperbolic approximation \eqref{eq:spatial_odd_hyperbolic}
to the KdV equation
\begin{equation}
\label{eq:kdv}
  \partial_t u(t, x) + \partial_x \frac{u(t,x)^2}{2} + \partial_x^3 u(t, x) = 0
\end{equation}
numerically. Thus, we discretize
the KdV equation as
\begin{equation}
  \partial_t \vec{u}
  + \frac{1}{3} \left(
    \vec{u} D_1 \vec{u} + D_1 \vec{u}^2
  \right)
  + D_+ D_0 D_- \vec{u}
  =
  \vec{0},
\end{equation}
where $D_\pm$ are again periodic upwind SBP operators and $D_0 = (D_+ + D_-)/2$
the corresponding central SBP operator. We use ARS(4,4,3) \cite{ascher1997implicit}
to integrate the resulting system in time, where we treat the nonlinear term
explicitly and the linear term implicitly.
Similarly, we discretize the hyperbolic approximation
\eqref{eq:spatial_odd_hyperbolic} as
\begin{equation}
  \partial_t
  \begin{pmatrix}
    \vec{q_0} \\
    \vec{q_1} \\
    \vec{q_2}
  \end{pmatrix}
  + \begin{pmatrix}
      \frac{1}{3} (\vec{u} D_1 \vec{u} + D_1 \vec{u}^2) \\
      \vec{0} \\
      \vec{0} \\
  \end{pmatrix}
  + \begin{pmatrix}
    D_+ \vec{q_2} \\
    \tau^{-1} (-D_0 \vec{q_1} + \vec{q_2}) \\
    \tau^{-1} (D_- \vec{q_0} - \vec{q_1})
  \end{pmatrix}
  =
  \begin{pmatrix}
    \vec{0} \\
    \vec{0} \\
    \vec{0}
  \end{pmatrix},
\end{equation}
where we apply the explicit part of the time integration scheme to the first term
and the implicit part to the second term. We initialize the hyperbolic approximation
with
\begin{equation}
  \vec{q_0} = \vec{u},
  \quad
  \vec{q_1} = D_- \vec{q_0} - \mu \vec{q_0},
  \quad
  \vec{q_2} = D_0 \vec{q_1}.
\end{equation}
The semidiscretization conserves discrete analogs of $\int q_0 \dif x$ and the energy
$\int (q_0^2 + \tau q_1^2 + \tau q_2^2) \dif x / 2$.
The asymptotic- and structure-preserving
properties of these discretizations are analyzed in \cite{biswas2025traveling}.

\begin{figure}[htb]
  \includegraphics[width=\textwidth]{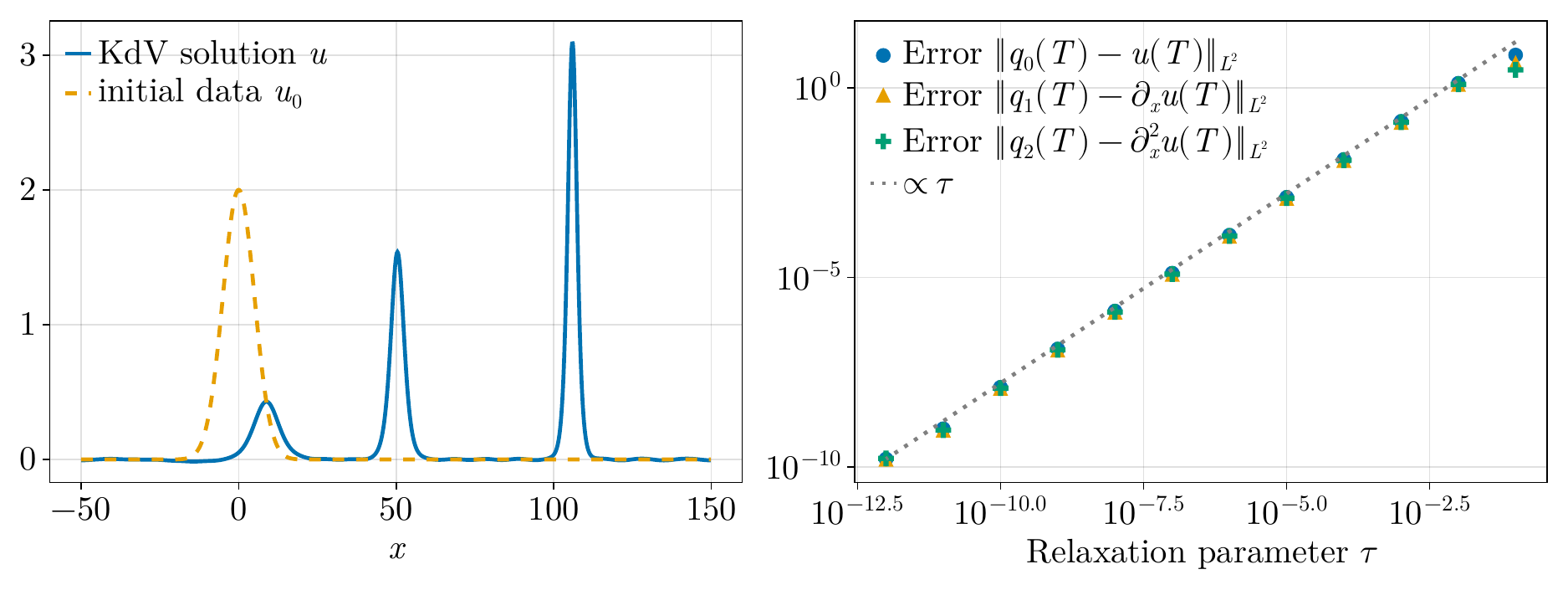}
  \caption{Convergence of the hyperbolic approximation
           \eqref{eq:spatial_odd_hyperbolic} with $\sigma_0 = 1$ and $\mu = 0$
           to the KdV equation \eqref{eq:kdv}. The left plot shows the numerical
           solution at the final time $T = 100$ and the initial condition.
           The right plot shows the convergence of the discrete $L^2$ error at the final time $T$
           as a function of the hyperbolic relaxation parameter $\tau$.}
  \label{fig:korteweg_de_vries_convergence}
\end{figure}

To test convergence of the hyperbolic approximation, we use the initial data
\begin{equation}
  u(0, x) = 2 \exp\bigl(-0.02 x^2\bigr)
\end{equation}
in the domain $\Omega = [-50, 150]$ with periodic boundary conditions.
The numerical solution of the KdV equation as well as the discrete $L^2$
errors are shown in Figure~\ref{fig:korteweg_de_vries_convergence}
for $2^{10}$ grid points, seventh-order accurate upwind operators, and
time step size $\Delta t = 0.05$.
As expected, the Gaussian initial condition splits into traveling waves at
the final time $t = 100$. Moreover, the hyperbolic approximation converges
to the KdV solution with the expected order $O(\tau)$ predicted by
Theorem~\ref{thm:convergence_spatial_odd}.
Again, the derivative approximations converge with the same order instead of
a reduced order as one might expect from the analysis in Theorem~\ref{thm:convergence_spatial_odd}.

\subsubsection{Korteweg-de Vries-Burgers equation}

Next, we consider the Korteweg-de Vries-Burgers equation
\begin{equation}
\label{eq:kdv_burgers}
  \partial_t u(t, x)
  + \partial_x \frac{u(t,x)^2}{2}
  - \mu \partial_{xx} u(t,x)
  + \partial_x^3 u(t, x) = 0
\end{equation}
with $\mu > 0$. Similar to the KdV equation, we discretize the KdV-Burgers equation as
\begin{equation}
  \partial_t \vec{u}
  + \frac{1}{3} \left(
    \vec{u} D_1 \vec{u} + D_1 \vec{u}^2
  \right)
  - \mu D_+ D_- \vec{u}
  + D_+ D_0 D_- \vec{u}
  =
  \vec{0},
\end{equation}
and integrate the nonlinear part explicitly while we treat the linear stiff part implicitly.
The discretization of the hyperbolic approximation \eqref{eq:spatial_odd_hyperbolic} reads
\begin{equation}
  \partial_t
  \begin{pmatrix}
    \vec{q_0} \\
    \vec{q_1} \\
    \vec{q_2}
  \end{pmatrix}
  + \begin{pmatrix}
      \frac{1}{3} (\vec{u} D_1 \vec{u} + D_1 \vec{u}^2) \\
      \vec{0} \\
      \vec{0} \\
  \end{pmatrix}
  + \begin{pmatrix}
    D_+ \vec{q_2} \\
    \tau^{-1} (-D_0 \vec{q_1} + \vec{q_2} + \mu \vec{q_1}) \\
    \tau^{-1} (D_- \vec{q_0} - \vec{q_1})
  \end{pmatrix}
  =
  \begin{pmatrix}
    \vec{0} \\
    \vec{0} \\
    \vec{0}
  \end{pmatrix}.
\end{equation}
The spatial semidiscretizations are energy-stable with dissipation term
mimicking the continuous case $-\mu \int (\partial_x u)^2 \dif x$.
Asymptotic-preserving properties can be analyzed similarly to the KdV
case in \cite{biswas2025traveling}.

\begin{figure}[htb]
  \includegraphics[width=\textwidth]{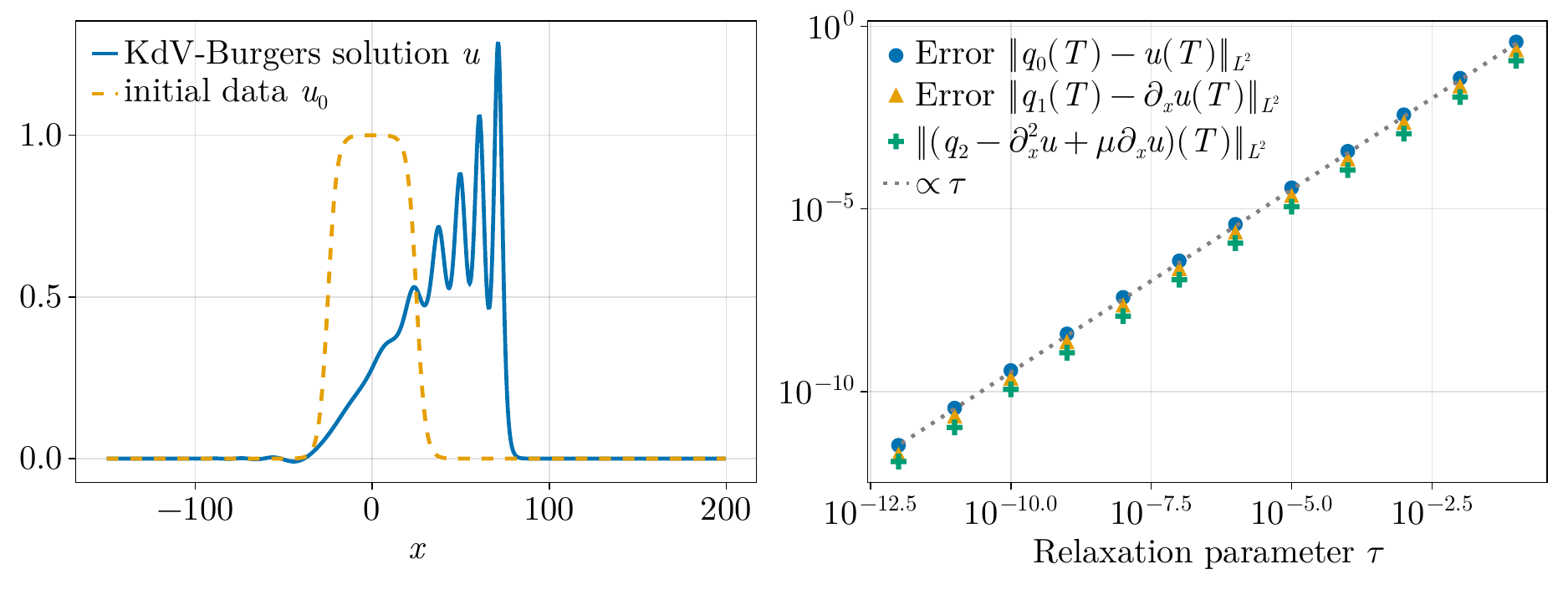}
  \caption{Convergence of the hyperbolic approximation
           \eqref{eq:spatial_odd_hyperbolic} with $\sigma_0 = 1$ and $\mu = 0.1$
           to the KdV-Burgers equation \eqref{eq:kdv_burgers}.
           The left plot shows the numerical solution at the final time
           $T = 100$ and the initial condition.
           The right plot shows the convergence of the discrete $L^2$ error at the final time $T$
           as a function of the hyperbolic relaxation parameter $\tau$.}
  \label{fig:korteweg_de_vries_burgers_convergence}
\end{figure}

We use the setup of \cite[Example~3.3]{xu2004local} with initial condition
\begin{equation}
  u(0, x) = \frac{1}{2} \left( 1 - \tanh\biggl( \frac{|x| - 25}{5} \biggr) \right)
\end{equation}
in the domain $\Omega = [-150, 200]$ with periodic boundary conditions.
The numerical solution of the KdV-Burgers equation as well as the discrete $L^2$
errors are shown in Figure~\ref{fig:korteweg_de_vries_burgers_convergence}
for $2^{10}$ grid points, seventh-order accurate upwind operators,
time step size $\Delta t = 0.1$, and dissipation parameter $\mu = 0.1$.
As expected, the hyperbolic approximation converges to the KdV-Burgers solution
with the expected order $O(\tau)$ predicted by Theorem~\ref{thm:convergence_spatial_odd}.
Again, the derivative approximations converge with the same order.

\subsubsection{Gardner equation}

We consider the Gardner equation
\begin{equation}
\label{eq:gardner}
  \partial_t u(t, x)
  + \sigma \partial_x \frac{u(t,x)^2}{2}
  + \partial_x \frac{u(t,x)^3}{3}
  + \partial_x^3 u(t, x) = 0
\end{equation}
and discretize it as
\begin{equation}
  \partial_t \vec{u}
  + \sigma \frac{1}{3} \left(
    \vec{u} D_1 \vec{u} + D_1 \vec{u}^2
  \right)
  + \frac{1}{6} \left( \vec{u}^2 D_1 \vec{u} + \vec{u} D_1 \vec{u}^2 + D_1 \vec{u}^3 \right)
  + D_+ D_0 D_- \vec{u}
  =
  \vec{0},
\end{equation}
and the hyperbolic approximation correspondingly. The split form of the cubic term
conserves the linear and quadratic invariants $\int u \dif x$ and $\int u^2 \dif x$
\cite{lefloch2021kinetic}.
Asymptotic-preserving properties can be analyzed similarly to the KdV
case in \cite{biswas2025traveling}.

\begin{figure}[htb]
  \includegraphics[width=\textwidth]{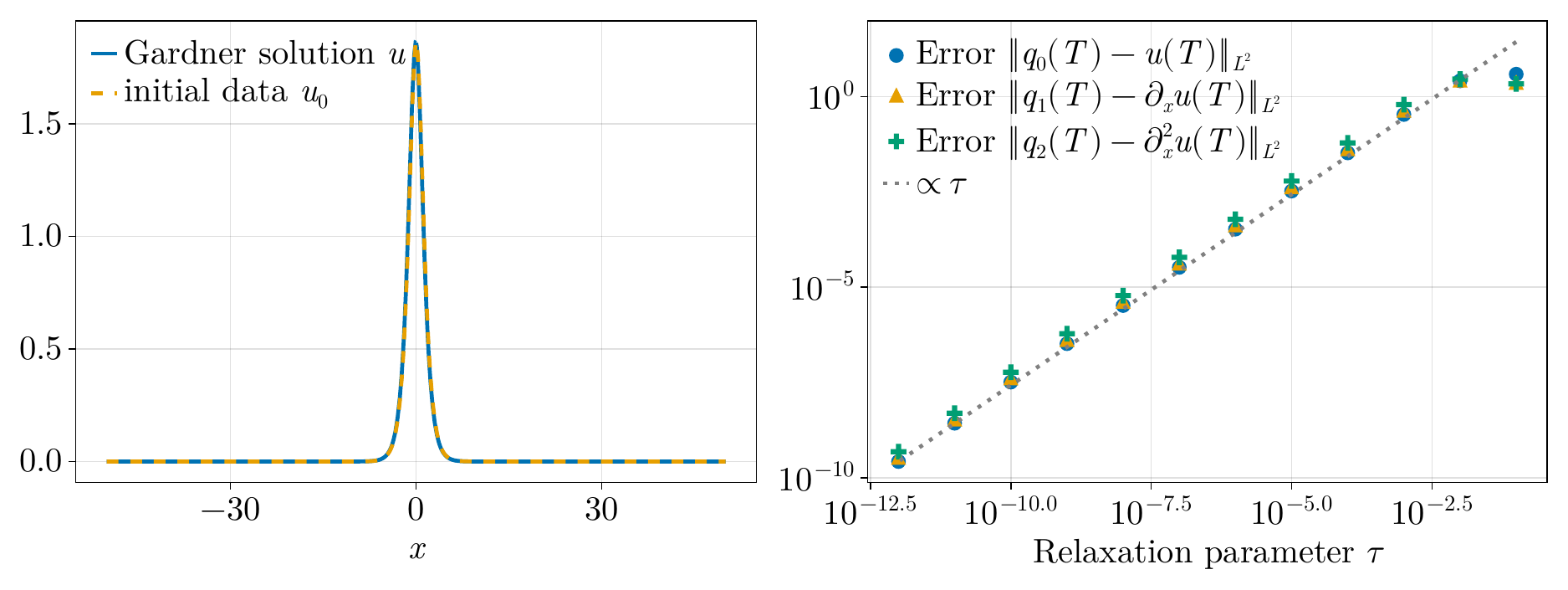}
  \caption{Convergence of the hyperbolic approximation
           to the Gardner equation \eqref{eq:gardner} with $\sigma = 1$.
           The left plot shows the numerical solution at the final time
           $T = 83.\overline{3}$ and the initial condition.
           The right plot shows the convergence of the discrete $L^2$ error at the final time $T$
           as a function of the hyperbolic relaxation parameter $\tau$.}
  \label{fig:gardner_convergence}
\end{figure}

We initialize the simulation with the value of the solitary wave solution
\begin{equation}
\begin{gathered}
  u(t,x) = \frac{A_1}{A_2 + \cosh\bigl(\sqrt{c} (x - c t) / 2\bigr)^2},
  \\
  A_1 = \frac{3 c}{\sqrt{\sigma^2 + 6 c}},
  \qquad
  A_2 = \frac{1}{2} \left( \frac{\sigma}{\sqrt{\sigma^2 + 6 c}} - 1 \right)
\end{gathered}
\end{equation}
in the domain $\Omega = [-50, 50]$ with periodic boundary conditions,
where $c = 1.2$ is the wave speed and $\sigma = 1$. The numerical solution
of the Gardner equation as well as the discrete $L^2$ errors are shown in
Figure~\ref{fig:gardner_convergence} for local discontinuous Galerkin (LDG) methods
using polynomials of degree $p = 8$ in 50 elements with time step size $\Delta t = 0.01$.
As expected, the hyperbolic approximation converges to the Gardner solution
with the expected order $O(\tau)$ predicted by Theorem~\ref{thm:convergence_spatial_odd}.
Again, the derivative approximations converge with the same order.

\subsection{Kawahara equation}
\label{sec:kawahara}

The Kawahara equation \cite{kawahara1972oscillatory} is
\begin{equation}
\label{eq:kawahara}
  \partial_t u + \partial_x f(u) + \partial_x^3 u - \partial_x^5 u = 0
\end{equation}
with $f(u) = u^2 / 2$.
Compared to \eqref{eq:spatial_odd_limit}, we have an additional dispersive
odd-order term $\partial_x^3 u$. Following the idea of
\cite[Section~4.3]{ketcheson2025approximation}, we could base the hyperbolic
approximation on the leading order $m = 5$ with $\sigma_0 = -1$
and approximate the third-order derivative $\partial_x^3 u$ by
$\partial_x q_2$, resulting in
\begin{equation}
\begin{aligned}
  \partial_t q_0 + \partial_x f(q_0) + \partial_x q_2 - \partial_x q_4 &= 0, \\
  \tau \partial_t q_1 + \partial_x q_3 &= q_4, \\
  \tau \partial_t q_2 - \partial_x q_2 &= -q_3, \\
  \tau \partial_t q_3 + \partial_x q_1 &= q_2, \\
  \tau \partial_t q_4 - \partial_x q_0 &= -q_1.
\end{aligned}
\end{equation}
For this system the functional $\eta$ defined by
\begin{equation}
 \eta(q)= \int_\Omega \left( \frac{1}{2} q_0^2 + \frac{1}{2} \tau \sum_{i=1}^4 q_i^2 \right) \dif x
\end{equation}
satisfies
\[
 \frac{\dif \eta(q)}{\dif t} = \int_\Omega - q_0 \partial_x q_2 \dif x,
\]
i.e., we cannot control its time evolution and should not call it an energy of this system. Similarly, we cannot set up a working relative energy framework for this equation
A seemingly straightforward alternative is to approximate $\partial_x^3 u$ by $q_3$
instead of $\partial_x q_2$ and obtain
\begin{equation}
\begin{aligned}
  \partial_t q_0 + \partial_x f(q_0) + q_3 - \partial_x q_4 &= 0, \\
  \tau \partial_t q_1 + \partial_x q_3 &= q_4, \\
  \tau \partial_t q_2 - \partial_x q_2 &= -q_3, \\
  \tau \partial_t q_3 + \partial_x q_1 &= q_2, \\
  \tau \partial_t q_4 - \partial_x q_0 &= -q_1.
\end{aligned}
\end{equation}
For this system
\[
 \frac{\dif \eta(q)}{\dif t} = \int_\Omega - q_0  q_3 \dif x,
\]
which means we have the same type of problem.
Instead, we propose the hyperbolic approximation
\begin{equation}
\label{eq:kawahara_hyperbolic}
\begin{aligned}
  \partial_t q_0 + \partial_x f(q_0) - \partial_x q_4 &= 0, \\
  \tau \partial_t q_1 - \partial_x q_1 + \partial_x q_3 &= q_4 , \\
  \tau \partial_t q_2 - \partial_x q_2 &= -q_3, \\
  \tau \partial_t q_3 + \partial_x q_1 &= q_2, \\
  \tau \partial_t q_4 - \partial_x q_0 &= -q_1
\end{aligned}
\end{equation}
whose solutions satisfy
\[
 \frac{\dif \eta(q)}{\dif t} = 0
\]
so that we can think of $\eta$ as an energy.
This system is hyperbolic with flux Jacobian
\begin{equation}
  \begin{pmatrix}
    q_0 & 0 & 0 & 0 & -1 \\
    0 & -\tau^{-1} & 0 & \tau^{-1} & 0 \\
    0 & 0 & -\tau^{-1} & 0 & 0 \\
    0 & \tau^{-1} & 0 & 0 & 0 \\
    -\tau^{-1} & 0 & 0 & 0 & 0
  \end{pmatrix},
\end{equation}
which is diagonalizable with the eigenvalues and eigenvectors given by
\begin{equation}
\begin{gathered}
  \lambda_1 = -\frac{1 + \sqrt{5}}{2 \tau}, \;
  \lambda_2 = \frac{-1}{\tau}, \;
  \lambda_3 = \frac{-1 + \sqrt{5}}{2 \tau}, \;
  \lambda_4 = \frac{q_0 - \sqrt{q_0^2 + 4 / \tau}}{2}, \;
  \lambda_5 = \frac{q_0 + \sqrt{q_0^2 + 4 / \tau}}{2}, \\
  v_1 = \begin{pmatrix} 0 \\ \tau \lambda_1 \\ 0 \\ 1 \\ 0 \end{pmatrix}, \quad
  v_2 = \begin{pmatrix} 0 \\ 0 \\ 1 \\ 0 \\ 0 \end{pmatrix}, \quad
  v_3 = \begin{pmatrix} 0 \\ \tau \lambda_3 \\ 0 \\ 1 \\ 0 \end{pmatrix}, \quad
  v_4 = \begin{pmatrix} -\tau \lambda_4 \\ 0 \\ 0 \\ 0 \\ 1 \end{pmatrix}, \quad
  v_5 = \begin{pmatrix} -\tau \lambda_5 \\ 0 \\ 0 \\ 0 \\ 1 \end{pmatrix}.
\end{gathered}
\end{equation}

To study the convergence of the hyperbolic approximation \eqref{eq:kawahara_hyperbolic},
we assume to be given a smooth solution $u$ to the Kawahara equation \eqref{eq:kawahara}.
Let $\bar q$ solve
\begin{equation}
\label{eq:kawahara_res}
\begin{aligned}
  \partial_t \bar q_0 + \partial_x f(\bar  q_0) - \partial_x \bar q_4 &= \tau R, \\
  \tau \partial_t \bar q_1 - \partial_x \bar q_1 + \partial_x \bar q_3 &=\bar  q_4 , \\
  \tau \partial_t \bar q_2 - \partial_x\bar  q_2 &= -\bar q_3, \\
  \tau \partial_t \bar q_3 + \partial_x \bar q_1 &= \bar q_2, \\
  \tau \partial_t \bar q_4 - \partial_x \bar q_0 &= -\bar q_1,
\end{aligned}
\end{equation}
with residual $R$ depending on $u$ and its derivatives and being uniformly bounded
for $\tau \to 0$. This can be achieved by setting $\bar q_2 =\partial_{xx} u$.
Then we compute $\bar q_3$ from \eqref{eq:kawahara_res}$_3$, i.e.,
\begin{equation}
  \bar q_3 = \partial_{xxx} u - \tau \partial_{txx} u.
\end{equation}
Let us also define primitive functions
$Q_3 = \partial_{xx} u - \tau \partial_{tx} u$ and
$P_3 = \partial_{x} u - \tau \partial_{t} u$.
Then we define $\bar q_1$ from \eqref{eq:kawahara_res}$_4$, i.e.
$\bar q_1 = \partial_x u + \tau \partial_t Q_3$ and its primitive
$Q_1 = u + \tau \partial_t P_3$.
Then we define $\bar q_4$ from \eqref{eq:kawahara_res}$_2$, i.e.,
\begin{equation}
  \bar q_4 = \partial_{xxxx} u - \tau \partial_{txxx} u
 - \partial_{xx} u - \tau \partial_{tx} Q_3
\end{equation}
and its primitive
$ Q_4=\partial_{xxx} u - \tau \partial_{txx} u
 - \partial_{x} u - \tau \partial_{t} Q_3$.
Finally, we define $\bar q_0$ from \eqref{eq:kawahara_res}$_5$, i.e.,
$\bar q_0= u + \tau \partial_t P_3- \tau \partial_t Q_4$.
In this way we get $\bar q_j = \partial_x^j u + \tau R_j$,
with $R_j$ depending on $u$ and its derivatives for $j \in \{0,\dots,3\}$
and $\bar q_4 = - \partial_{xx} u + \partial_x^4 u + \tau R_4$.

Thus, the residual $R$ in \eqref{eq:kawahara_res} satisfies
\begin{multline}
 |R|= |\partial_t ( \partial_t P_3 - \partial_t Q_4) - \tau \partial_x R_4
 + \partial_x (f(\bar q_0) - f(u))|
 \\
 \leq
 |\partial_t ( \partial_t P_3 - \partial_t Q_4) - \tau \partial_x R_4|
 + \|f''\|_\infty |\partial_t P_3 - \partial_t Q_4|  \, |\partial_x \bar q_0|
 + |f'(u)| |\partial_{xt} P_3 - \partial_{tx} Q_4|.
\end{multline}

We define the relative energy
\begin{equation}
  \eta(q, \bar q)
  =
  \int \left(
    \frac12 (q_0 - \bar q_0)^2 + \frac{\tau}{2} \sum_{j=1}^3 (q_j - \bar q_j)^2
  \right) \dif x.
\end{equation}
Its time evolution can be computed setting $f(u) = u^2 / 2$ as
\begin{equation}
\begin{aligned}
  &\qquad
  \frac{\dif}{\dif t} \eta(q,\bar q)
  =
  \int \left( (q_0- \bar q_0) \partial_t (q_0- \bar q_0) + \sum_{j=1}^{3} \tau (q_j- \bar q_j) \partial_t (q_j- \bar q_j) \right) \dif x
  \\
  &= -\int \biggl(-(q_0- \bar q_0) \partial_x(f(q_0)-f(\bar q_0)) +
   (q_0- \bar q_0)\partial_x (q_4- \bar q_4)+  (q_0- \bar q_0)\tau R\\
   & \qquad + (q_1- \bar q_1)\partial_x(q_1 - \bar q_1)
   -  (q_1- \bar q_1)\partial_x(q_3 - \bar q_3) + (q_1- \bar q_1)(q_4 - \bar q_4) +  (q_2- \bar q_2)\partial_x(q_2 - \bar q_2)\\
   &\qquad - (q_2- \bar q_2)(q_3 - \bar q_3)-  (q_3- \bar q_3)\partial_x(q_1 - \bar q_1)
   +  (q_3- \bar q_3)(q_2 - \bar q_2)+ (q_4- \bar q_4)\partial_x(q_0 - \bar q _0)\\
   &\qquad- (q_4- \bar q_4)(q_1 - \bar q_1) \biggr)\dif x\\
  &= \int \biggl(- (q_0- \bar q_0) \partial_x(f(q_0)-f(\bar q_0)) +  (q_0- \bar q_0)\tau R \biggr) \dif x.
\end{aligned}
\end{equation}
From here we can proceed as before and obtain
\begin{equation}
 \| u - q_0 \|_{L^\infty(0,T, L^2(\Omega))}+ \sum_{j=1}^{3} \sqrt{\tau} \| \partial_x^j u - q_j \|_{L^\infty(0,T, L^2(\Omega))} = \mathcal{O}(\tau).
\end{equation}
Thus, we obtain
\begin{theorem}
\label{thm:convergence_kawahara}
  Let $T>0$ and
  let $u \in W^{5,\infty}((0,T) \times \Omega)$ be a solution to
  \eqref{eq:kawahara} with initial data $u_0 \in W^{5,\infty}(\Omega)$.
  Let, for each $\tau>0$, $q$ be an entropy solution to \eqref{eq:kawahara_hyperbolic}
  with $q|_{t=0}= (u_0, \partial_x u_0 ,\dots, \partial_x^{4} u_0 )$. Then
\begin{equation}
 \| u - q_0 \|_{L^\infty(0,T, L^2(\Omega))}+ \sum_{j=1}^{3} \sqrt{\tau} \| \partial_x^j u - q_j \|_{L^\infty(0,T, L^2(\Omega))} = \mathcal{O}(\tau).
\end{equation}
\end{theorem}

\begin{remark}
\label{rem:convergence_kawahara}
  It should be noted that $f(u)=u^2/2$ is not needed for the theorem to hold.
  What is needed is that there exists $C>0$ such that
  \begin{equation*}
    |f(q_0) - f(\bar q_0) - f'(\bar q_0)(q_0 - \bar q_0)|
    \leq
    C |q_0 - q_0|^2 \quad \text{ uniformly in } (0,T) \times \Omega.
  \end{equation*}
  This means theorem also holds for generalized Kawahara equations with $f'' \in L^{\infty}(\mathbb{R})$  or with $f \in C^2(\mathbb{R})$ provided there exists some compact $K \subset \mathbb{R}$ such that $\bar q_0$ (as constructed from $u$) and $q_0$ (for all values of $\tau$) only take values in this set.
\end{remark}

\subsubsection{Numerical demonstration}

We discretize the Kawahara equation \eqref{eq:kawahara} as
\begin{equation}
  \partial_t \vec{u}
  + \frac{1}{3} \left( \vec{u} D_1 \vec{u} + D_1 \vec{u}^2 \right)
  + D_+ D_0 D_- \vec{u}
  - D_+ D_+ D_0 D_- D_- \vec{u}
  =
  \vec{0},
\end{equation}
and the hyperbolic approximation \eqref{eq:kawahara_hyperbolic} as
\begin{equation}
  \partial_t
  \begin{pmatrix}
    \vec{q_0} \\
    \vec{q_1} \\
    \vec{q_2} \\
    \vec{q_3} \\
    \vec{q_4}
  \end{pmatrix}
  + \begin{pmatrix}
      \frac{1}{3} (\vec{q_0} D_1 \vec{q_0} + D_1 \vec{q_0}^{\!\!2}) \\
      \vec{0} \\
      \vec{0} \\
      \vec{0} \\
      \vec{0} \\
  \end{pmatrix}
  + \begin{pmatrix}
    D_+ \vec{q_4} \\
    \tau^{-1} (D_+ \vec{q_3} - D_0 \vec{q_1} - \vec{q_4}) \\
    \tau^{-1} (-D_0 \vec{q_2} + \vec{q_3}) \\
    \tau^{-1} (D_- \vec{q_1} - \vec{q_2}) \\
    \tau^{-1} (-D_- \vec{q_0} + \vec{q_1})
  \end{pmatrix}
  =
  \begin{pmatrix}
    \vec{0} \\
    \vec{0} \\
    \vec{0} \\
    \vec{0} \\
    \vec{0}
  \end{pmatrix}.
\end{equation}
We treat the nonlinear term explicitly and the linear term implicitly
using ARS(4,4,3) \cite{ascher1997implicit}.
The spatial semidiscretizations are energy-conservative.
Asymptotic-preserving properties can be analyzed similarly to the KdV
case in \cite{biswas2025traveling}.

\begin{figure}[htb]
  \includegraphics[width=\textwidth]{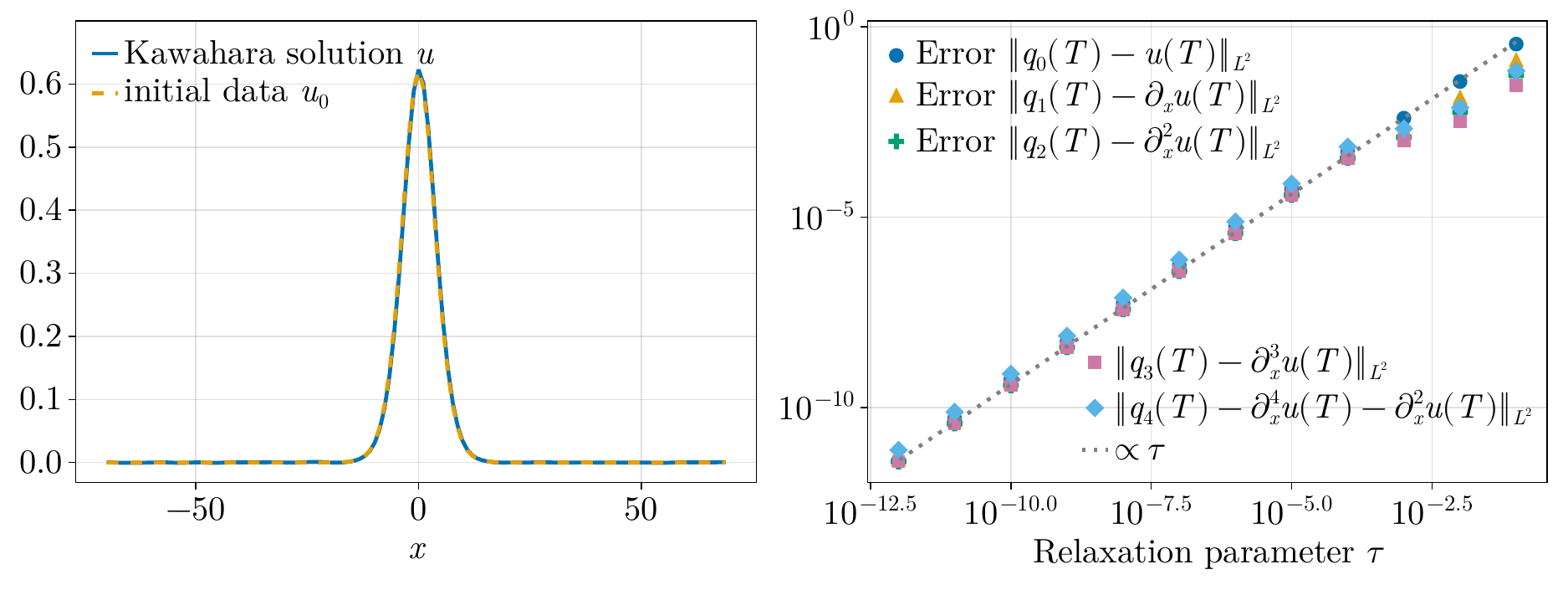}
  \caption{Convergence of the hyperbolic approximation
           \eqref{eq:kawahara_hyperbolic}
           to the Kawahara equation \eqref{eq:kawahara}.
           The left plot shows the numerical solution at the final time
           $T = 657.\overline{2}$ and the initial condition.
           The right plot shows the convergence of the discrete $L^2$ error at the final time $T$
           as a function of the hyperbolic relaxation parameter $\tau$.}
  \label{fig:kawahara_convergence}
\end{figure}

We use the setup of \cite[Example~3.4]{xu2004local} with initial condition
given by the exact solitary wave solution
\begin{equation}
  u(t, x) = \frac{105}{169} \operatorname{sech}\biggl(\frac{x - \frac{36}{169} t}{2 \sqrt{13}}\biggr)^4.
\end{equation}
We use periodic boundary conditions in the domain $[-70, 70]$
and discretize the problem with $2^{7}$ grid points, third-order accurate
upwind operators, and time step size $\Delta t = 0.1$.
The numerical solution of the Kawahara equation after a full traversal of the domain
as well as the discrete $L^2$ errors are shown in Figure~\ref{fig:kawahara_convergence}.
As expected, the hyperbolic approximation converges to the Kawahara solution
with the expected order $O(\tau)$ predicted by Theorem~\ref{thm:convergence_kawahara}.
Again, the derivative approximations converge with the same order.

\subsubsection{Generalized Kawahara equation}
\label{sec:generalized_kawahara}

Next, we consider the generalized Kawahara equation
\begin{equation}
\label{eq:generalized_kawahara}
  \partial_t u + \partial_x f(u) + \partial_x^3 u - \partial_x^5 u = 0
\end{equation}
with $f(u) = \sigma u^2 / 2 + u^3 / 3$ as in \cite[Example~3.5]{xu2004local}.
We discretize it as
\begin{equation}
  \partial_t \vec{u}
  + \frac{\sigma}{3} \left( \vec{u} D_1 \vec{u} + D_1 \vec{u}^2 \right)
  + \frac{1}{6} \left( \vec{u}^2 D_1 \vec{u} + \vec{u} D_1 \vec{u}^2 + D_1 \vec{u}^3 \right)
  + D_+ D_0 D_- \vec{u}
  - D_+ D_+ D_0 D_- D_- \vec{u}
  =
  \vec{0},
\end{equation}
and the hyperbolic approximation correspondingly. The split form of the cubic term
conserves the linear and quadratic invariants $\int u \dif x$ and $\int u^2 \dif x$
\cite{lefloch2021kinetic}.
Asymptotic-preserving properties can be analyzed similarly to the KdV
case in \cite{biswas2025traveling}.

\begin{figure}[htb]
  \includegraphics[width=\textwidth]{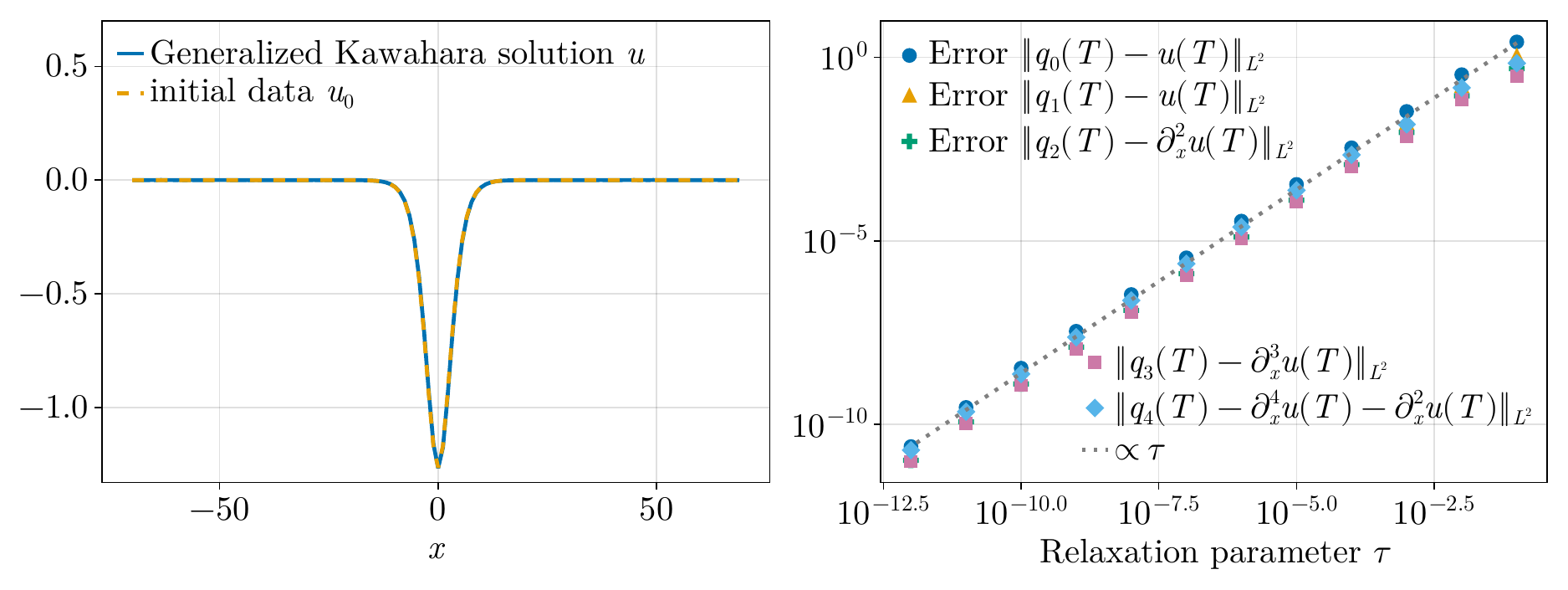}
  \caption{Convergence of the hyperbolic approximation
           to the generalized Kawahara equation \eqref{eq:generalized_kawahara}.
           The left plot shows the numerical solution at the final time
           $T = 715.\overline{90}$ and the initial condition.
           The right plot shows the convergence of the discrete $L^2$ error at the final time $T$
           as a function of the hyperbolic relaxation parameter $\tau$.}
  \label{fig:generalized_kawahara_convergence}
\end{figure}

We initialize the generalized Kawahara equation with the exact solitary wave
solution \cite[Example~3.5]{xu2004local}
\begin{equation}
  u(t, x) = -6 \sqrt{10} k^2 \operatorname{sech}\bigl( k (x -c t) \bigr)^2,
\end{equation}
where $k = \sqrt{1 / 20 + \sigma / (4 \sqrt{10})}$, $c = 4 k^2 (1 - 4 k^2)$,
and $\sigma = 2 / \sqrt{90}$.
We use periodic boundary conditions in the domain $[-70, 70]$
and discretize the problem with $2^{7}$ grid points, seventh-order accurate
upwind operators, and time step size $\Delta t = 0.1$.
The numerical solution of the generalized Kawahara equation after a
full traversal of the domain as well as the discrete $L^2$ errors are shown in
Figure~\ref{fig:generalized_kawahara_convergence}.
As expected, the hyperbolic approximation converges to the generalized Kawahara solution
with the expected order $O(\tau)$.
Again, the derivative approximations converge with the same order.

\subsection{Relative equilibrium structure}

The hyperbolic approximations of the energy-conservative PDEs discussed above
have a specialized structure that allows to obtain significantly improved
error growth rates in time for solitary wave solutions. We will describe this
so-called relative equilibrium structure in the following for the generalized
Kawahara equation \eqref{eq:generalized_kawahara} and its hyperbolic approximation.

The generalized Kawahara equation \eqref{eq:generalized_kawahara} is a Hamiltonian PDE
\begin{equation}
\label{eq:hamiltonian_pde}
  \partial_t u = \mathcal{J} \, \delta \mathcal{H}(u)
\end{equation}
with skew-symmetric operator $\mathcal{J} = -\partial_x$ and Hamiltonian
\begin{equation}
  \mathcal{H}(u) = \int \left( F(u) + \frac{1}{2} (\partial_x u)^2 - \frac{1}{2} (\partial_x^2 u)^2 \right) \dif x,
\end{equation}
where $F(u)$ is the primitive of $f(u) = \sigma u^2 / 2 + u^3 / 3$, i.e.,
\begin{equation}
  F(u) = \frac{\sigma}{6} u^3 + \frac{1}{12} u^4.
\end{equation}
Moreover, the generalized Kawahara equation has the invariant
\begin{equation}
  \mathcal{I}(u) = \int \frac{1}{2} u^2 \dif x,
\end{equation}
satisfying
\begin{equation}
  \mathcal{J} \, \delta \mathcal{I}(u) = -\partial_x u.
\end{equation}
This is the relative equilibrium structure discussed by
Dur{\'a}n and Sanz-Serna \cite{duran1998numerical,duran2000numerical}. This theory predicts
that the error of solitary waves grows quadratically in time for general methods but
only linearly in time for time integration schemes that conserve the Hamiltonian $\mathcal{H}$
or the invariant $\mathcal{I}$ \cite{frutos1997accuracy}. To make the results rigorous,
PDE stability estimates are required, which we will not discuss here. However, the general
theory has been applied successfully to several PDEs, e.g.,
\cite{araujo2001error,duran2002numerical,ranocha2021rate,ranocha2025structure}.

The hyperbolic approximation of the generalized Kawahara equation
\eqref{eq:generalized_kawahara} inherits the relative equilibrium structure. Indeed,
\begin{equation}
\label{eq:generalized_kawahara_hyperbolic}
\begin{aligned}
  \partial_t q_0 + \partial_x f(q_0) - \partial_x q_4 &= 0, \\
  \tau \partial_t q_1 - \partial_x q_1 + \partial_x q_3 &= q_4 , \\
  \tau \partial_t q_2 - \partial_x q_2 &= -q_3, \\
  \tau \partial_t q_3 + \partial_x q_1 &= q_2, \\
  \tau \partial_t q_4 - \partial_x q_0 &= -q_1,
\end{aligned}
\end{equation}
can be written as Hamiltonian PDE \eqref{eq:hamiltonian_pde} with Hamiltonian
\begin{equation}
  \mathcal{H}(q) = \int\left(
    F(q_0) - q_0 q_4 - \frac{q_1^2}{2} + q_1 q_3 - q_1 P(q_4) - \frac{q_2^2}{2} + q_2 P(q_3)
  \right) \dif x,
\end{equation}
where
\begin{equation}
  P(q_i)(t,x) = \int^x q_i(t, y) \dif y
\end{equation}
denotes the primitive function of $q_i$. Since
\begin{equation}
  \delta \mathcal{H}(q) = \begin{pmatrix}
    f(q_0) - q_4 \\
    -q_1 + q_3 - P(q_4) \\
    -q_2 + P(q_3) \\
    q_1 - P(q_2) \\
    -q_0 + P(q_1)
  \end{pmatrix},
\end{equation}
\eqref{eq:generalized_kawahara_hyperbolic} is obtained using the skew-symmetric
operator
\begin{equation}
  \mathcal{J} = \begin{pmatrix}
    -\partial_x  \\
    & -\tau^{-1}\partial_x \\
    & & -\tau^{-1}\partial_x \\
    & & & \tau^{-1}\partial_x \\
    & & & & \tau^{-1}\partial_x
  \end{pmatrix}.
\end{equation}
Moreover, the invariant
\begin{equation}
  \mathcal{I}(q) = \int \left(
    \frac{1}{2} q_0^2 + \frac{\tau}{2} \sum_{i=1}^4 q_i^2
  \right) \dif x
\end{equation}
satisfies
\begin{equation}
  \mathcal{J} \, \delta \mathcal{I}(q) = -\begin{pmatrix}
    \partial_x q_0 \\
    \partial_x q_1 \\
    \partial_x q_2 \\
    \partial_x q_3 \\
    \partial_x q_4
  \end{pmatrix}.
\end{equation}
Thus, the hyperbolic approximation shares the relative equilibrium structure,
and we can expect reduced error growth rates for solitary wave solutions.

\begin{figure}[htb]
  \includegraphics[width=\textwidth]{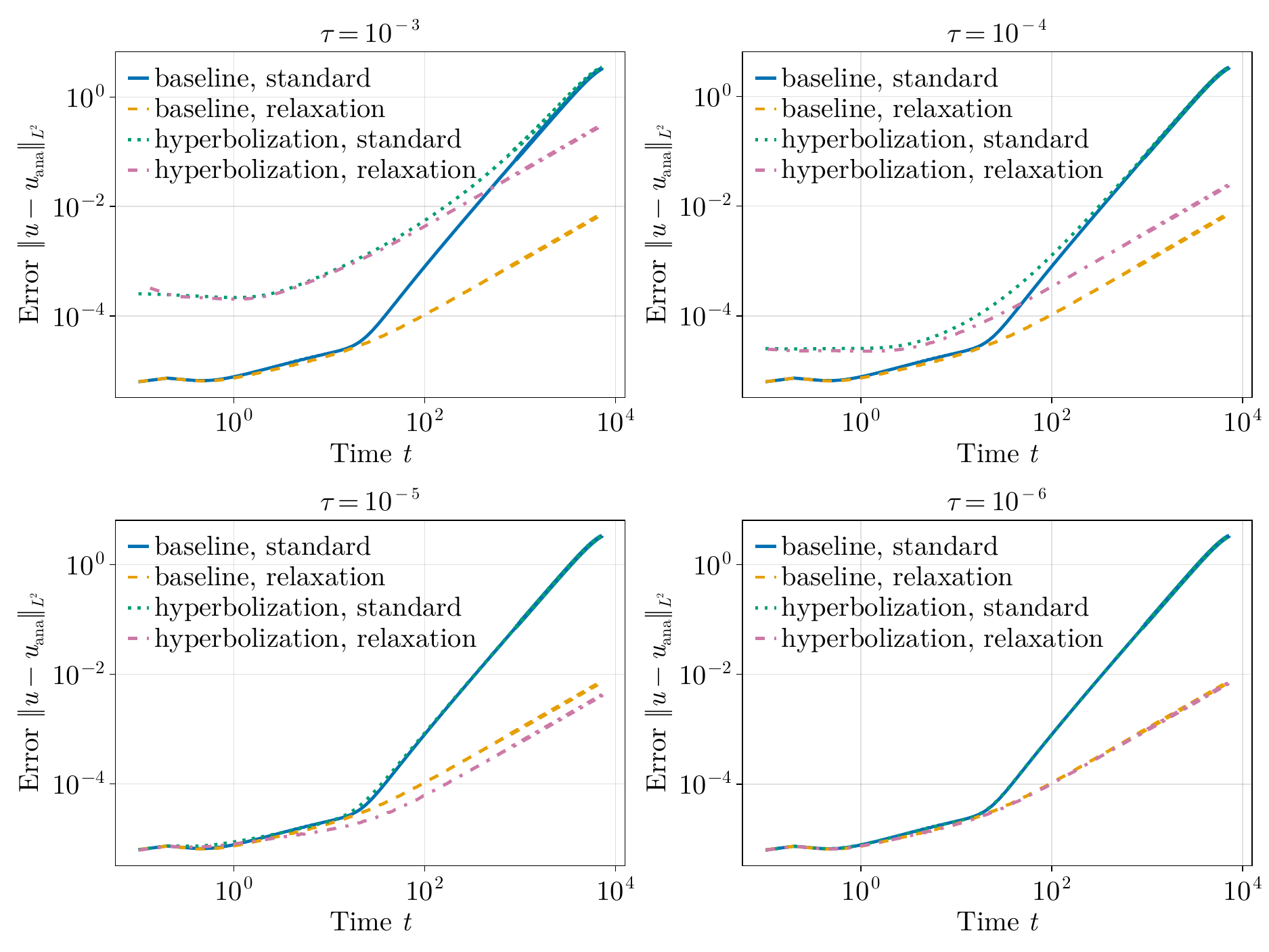}
  \caption{Error growth in time for numerical solutions of the
           generalized Kawahara equation \eqref{eq:generalized_kawahara}
           and its hyperbolic approximation for different values of the
           relaxation time $\tau$.}
  \label{fig:generalized_kawahara_error_growth}
\end{figure}

Numerical results obtained for the solitary wave setup described in
Section~\ref{sec:generalized_kawahara} are presented in
Figure~\ref{fig:generalized_kawahara_error_growth}. Here, we used
$2^9$ grid points in the domain $[-70, 70]$, seventh-order accurate
upwind operators, and time step size $\Delta t = 0.1$ to compute the
numerical solutions during ten full traversals of the domain.
We used the relaxation method
\cite{ketcheson2019relaxation,ranocha2020relaxation,ranocha2020general}
to conserve the quadratic invariant $\mathcal{I}$.
In all cases, we see that the relaxation method results in linear error
growth in time while the standard method results in quadratic error growth,
in accordance with the relative equilibrium theory of Dur{\'a}n and Sanz-Serna
\cite{duran1998numerical,duran2000numerical}.
For relatively big values $\tau \in \{10^{-3}, 10^{-4}\}$, the error of the
hyperbolic approximation is bigger than the error obtained by discretizing
the baseline system. Interestingly, the error of the hyperbolization with
$\tau = 10^{-5}$ and relaxation in time to conserve the quadratic invariant
is smaller than the error obtained by simulating the original PDE. For smaller
values of $\tau$, the error plots of the baseline PDE and the hyperbolization
become visually indistinguishable due to the convergence of the hyperbolic
approximation.

\subsection{Even leading order \texorpdfstring{$m$}{m}}

Next, we study PDEs of the form
\begin{equation}
\label{eq:spatial_even_limit}
  \partial_t u(t, x)
  + \partial_x f(u(t,x))
  + \sigma_0 \partial_x^m u(t, x) = 0
\end{equation}
with even $m \in \N$ and $\sigma_0 = (-1)^{m / 2}$ (as required for stability).
We generalize the hyperbolic approximations \eqref{eq:ketcheson_biswas_hyperbolic}
of \cite{ketcheson2025approximation} to
\begin{equation}
\begin{aligned}
\label{eq:spatial_even_hyperbolic}
  \partial_t q_{0}
  + \partial_x f(q_0)
  + \sigma_0 \partial_x q_{m-1}
  &= 0,
  \\
  \tau \partial_t q_{j}
  + \sigma_0 (-1)^{j} \partial_x q_{m-j-1}
  &= \sigma_0 (-1)^{j} q_{m-j},
  & j &\in \{1, \dots, \lceil m / 2 \rceil - 1\},
  \\
  \tau \partial_t q_{j}
  + \sigma_0 (-1)^{j+m-1} \partial_x q_{m-j-1}
  &= \sigma_0 (-1)^{j+m-1} q_{m-j},
  & j &\in \{\lceil m / 2 \rceil, \dots, m-1\}.
\end{aligned}
\end{equation}

To perform the convergence analysis using the relative energy method,
we assume to be given a smooth solution $u$ to \eqref{eq:spatial_even_limit}.
We wish to define $\bar q$ such that it satisfies
\begin{equation}
\label{eq:spatial_even_res}
\begin{aligned}
  \partial_t \bar q_{0}
  + \partial_x f(\bar q_{0})
  + \sigma_0 \partial_x \bar q_{m-1}
  &= \tau R,
  \\
  \tau \partial_t \bar q_{j}
  + \sigma_0 (-1)^{j} \partial_x \bar q_{m-j-1}
  &= \sigma_0 (-1)^{j} \bar q_{m-j},
  & j &\in \{1, \dots, \lceil m / 2 \rceil - 1\},
  \\
  \tau \partial_t \bar q_{j}
  + \sigma_0 (-1)^{j+m-1} \partial_x \bar q_{m-j-1}
  &= \sigma_0 (-1)^{j+m-1} \bar q_{m-j},
  & j &\in \{\lceil m / 2 \rceil, \dots, m-1\},
\end{aligned}
\end{equation}
where the residual $R$ depends on $u$ and its derivatives and is uniformly
bounded for $ \tau \to 0$. The crucial point is again that there is only a residual
in the first equation.

We start by setting $\bar q_{\frac{m}{2}}= \partial_x^{\frac{m}{2}}$ and consider the evolution equation for $q_j$ with $j= \frac{m}{2}$ to determine $q_{\frac{m}{2}-1}$.
Then we do the following two steps successively for all $k \in \{1, \dots, \frac{m}{2}\}$:
\begin{enumerate}
  \item We use the evolution equation for $q_j$ with $j=\frac{m}{2}-k$
        to determine $q_{\frac{m}{2}+k}$, i.e.,
        \begin{equation}
          \tau \partial_t q_{\frac{m}{2}-k} + \sigma_0 (-1)^{\frac{m}{2}-k} \partial_x q_{\frac{m}{2}+k-1} = \sigma_0 (-1)^{\frac{m}{2}-k} q_{\frac{m}{2}+k}.
        \end{equation}
        Note that $q_{\frac{m}{2}-k}$ and $q_{\frac{m}{2}+k-1}$ have already been
        determined in the iteration.
  \item We use the evolution equation for $q_j$ with $j=\frac{m}{2}+k$
        to determine $q_{\frac{m}{2}-k-1}$, i.e.,
        \begin{equation}
          \tau \partial_t q_{\frac{m}{2}+k} + \sigma_0 (-1)^{\frac{m}{2}+k-1} \partial_x q_{\frac{m}{2}-1-k} = \sigma_0 (-1)^{\frac{m}{2}+k-1} q_{\frac{m}{2}-k}
        \end{equation}
        Note that $q_{\frac{m}{2}+k}$ and $q_{\frac{m}{2}-k}$ have already been
        determined in the iteration and each term in them contains at least one $x$-derivative.
\end{enumerate}

By this procedure, $\bar q$ satisfies each equation in \eqref{eq:spatial_even_hyperbolic},
except the first one, exactly.
In addition, if we choose the integration constants in step 2 reasonably, for $j \in \{0, \dots, m-1\}$, each $\bar q_j$ satisfies
\begin{equation}
  \bar q_j= \partial_x^j u + \tau R_j,
\end{equation}
where $R_j$ depends on $u$ and its derivatives and is uniformly bounded for $\tau \to 0$.
Similar to the discussion for odd $m$, if $j = \tfrac{m-1\pm k}{k}$ with $k$ odd
and positive, then $R_j$ only depends on derivatives of $u$ up to order $\frac{m}{2} + \frac{k-1}{2}$.
Indeed, the definition of $\bar q_0$ will contain derivatives of $u$ of degree $m-1$
so that the convergence proof requires $u \in W^{m,\infty}((0,T) \times \Omega)$.

We define the relative energy
\begin{equation}
  \eta(q, \bar q)
  :=
  \int_\Omega \left(
    \frac12 (q_0 - \bar q_0)^2 + \frac{\tau}{2} \sum_{i=1}^{m-1} (q_i - \bar q_i)^2
  \right)  \dif x
\end{equation}
and, based on \eqref{eq:spatial_even_res}, we compute
\begin{equation}
\begin{aligned}
  \frac{\dif}{\dif t} \eta(q,\bar q)
  &=
  \int \left( (q_0- \bar q_0) \partial_t (q_0- \bar q_0) + \sum_{j=1}^{m-1} \tau (q_j- \bar q_j) \partial_t (q_j- \bar q_j) \right) \dif x
  \\
  &= -\int \biggl(
    (q_0- \bar q_0) \partial_x \big(f(q_0)- f(\bar q_0)\big) +(q_0- \bar q_0)\tau R + \sigma_0 (q_0- \bar q_0) \partial_x (q_{m-1}- \bar q_{m-1})\\
    &\qquad\qquad
    + \sum_{j=1}^{\lceil m/2 \rceil - 1} \sigma_0 (-1)^j (q_j- \bar q_j)
      \bigl( \partial_x q_{m-j-1} - q_{m-j} - \partial_x \bar q_{m-j-1} + \bar q_{m-j} \bigr)
  \\
  &\qquad\qquad
    + \sum_{j=\lceil m / 2 \rceil}^{m - 1} \sigma_0 (-1)^{j+m-1} (q_{j}- \bar q_j)
      \bigl( \partial_x q_{m-j-1} - q_{m-j}-  \partial_x \bar q_{m-j-1} + \bar q_{m-j} \bigr)
  \biggr) \dif x
  \\
  &= -\int \biggl(
   (q_0- \bar q_0) \partial_x \big(f(q_0)- f(\bar q_0)\big) +(q_0- \bar q_0)\tau R\\
    &\qquad\qquad
    + \sigma_0 \sum_{j=0}^{\lceil m/2 \rceil - 1} (-1)^j  (q_j- \bar q_j) \partial_x (q_{m-j-1}-\bar q_{m-j-1})
    \\
    &\qquad\qquad
    - \sigma_0 \sum_{j=\lceil m/2 \rceil }^{m - 1} (-1)^j  (q_j- \bar q_j) \partial_x (q_{m-j-1}- q_{m-j-1})
    \\
    &\qquad\qquad
    - \sigma_0 \sum_{j=1}^{\lceil m/2 \rceil - 1} (-1)^j  (q_j- \bar q_j) (q_{m-j}-\bar q_{m-j})
    \\
    &\qquad\qquad
    + \sigma_0 \sum_{j=\lceil m/2 \rceil}^{m - 1} (-1)^j  (q_j- \bar q_j) (q_{m-j} -\bar q_{m-j})
  \biggr) \dif x
  \\
  &=
  -\sigma_0 (-1)^{m/2} \int q_{m/2}^2 \dif x
  + \tau^2 \int R^2 \dif x + \eta(q, \bar q)\\
  &
  \le  \tau^2 \int R^2 \dif x + \|\partial_x  \bar q_0 \|_{L^\infty} \eta(q, \bar q).
\end{aligned}
\end{equation}
We can conclude
\begin{equation}
  \| \eta(q, \bar q)\|_{L^\infty(0,T)} = \mathcal{O}(\tau^2)
\end{equation}
by using Gronwall's lemma since $u$ is independent of $\tau$
and $\eta(q, \bar q)|_{t=0}= \mathcal{O}(\tau^2)$.
This implies, in particular, that
\begin{equation}
 \| u - q_0 \|_{L^\infty(0,T, L^2(\Omega))}+ \sqrt{\tau}\sum_{j=1}^{m-1}\| \partial_x^j u - q_j \|_{L^\infty(0,T, L^2(\Omega))} = \mathcal{O}(\tau).
\end{equation}
\begin{theorem}
\label{thm:convergence_spatial_even}
  Let $T>0$ and $f \in W^{2,\infty}_{loc}(\mathbb{R})$ such that
  $f'' \in L^\infty(\mathbb{R})$.
  Let $u \in W^{m,\infty}((0,T) \times \Omega)$ be a solution to
  \eqref{eq:spatial_even_limit} with initial data $u_0 \in W^{m,\infty}(\Omega)$.
  Let, for each $\tau>0$, $q$ be an entropy solution to \eqref{eq:spatial_even_hyperbolic}
  with $q|_{t=0}= (u_0, \partial_x u_0 ,\dots, \partial_x^{m-1} u_0 )$. Then
  \begin{equation}
    \| u - q_0 \|_{L^\infty(0,T, L^2(\Omega))}
    + \sum_{j=1}^{m-1} \sqrt{\tau} \| \partial_x^j u - q_j \|_{L^\infty(0,T, L^2(\Omega))}
    = \mathcal{O}(\tau).
  \end{equation}
\end{theorem}

\begin{remark}
 Note that an analogous statement to Remark \ref{rem:convergence_mixed} also holds in this case.
\end{remark}

\subsubsection{Linear bi-harmonic equation}

We consider the linear bi-harmonic equation \cite[Example~5.3]{yan2002local}
\begin{equation}
\label{eq:linear_biharmonic}
  \partial_t u + \partial_x^4 u = 0
\end{equation}
with solution $u(t, x) = \exp(-t) \sin(x)$ in the domain $\Omega = [0, 2\pi]$
with periodic boundary conditions. We discretize the linear bi-harmonic equation as
\begin{equation}
  \partial_t \vec{u} + D_+ D_- D_+ D_- \vec{u} = \vec{0},
\end{equation}
and the hyperbolic approximation \eqref{eq:spatial_even_hyperbolic} as
\begin{equation}
  \partial_t
  \begin{pmatrix}
    \vec{q_0} \\
    \vec{q_1} \\
    \vec{q_2} \\
    \vec{q_3}
  \end{pmatrix}
  + \begin{pmatrix}
      D_+ \vec{q_3} \\
      \tau^{-1} (-D_- \vec{q_2} + \vec{q_3}) \\
      \tau^{-1} (-D_+ \vec{q_1} + \vec{q_2}) \\
      \tau^{-1} (D_- \vec{q_0} - \vec{q_1})
  \end{pmatrix}
  =
  \begin{pmatrix}
    0 \\
    0 \\
    0 \\
    0
  \end{pmatrix}.
\end{equation}
We treat all linear terms implicitly.
The semidiscretization conserves a discrete analog of $\int q_0 \dif x$ and dissipates
the energy $\int (q_0^2 + \tau q_1^2 + \tau q_2^2 + \tau q_3^2) \dif x / 2$.
Asymptotic-preserving properties can be analyzed similarly to the KdV
case in \cite{biswas2025traveling}.

\begin{figure}[htb]
  \includegraphics[width=\textwidth]{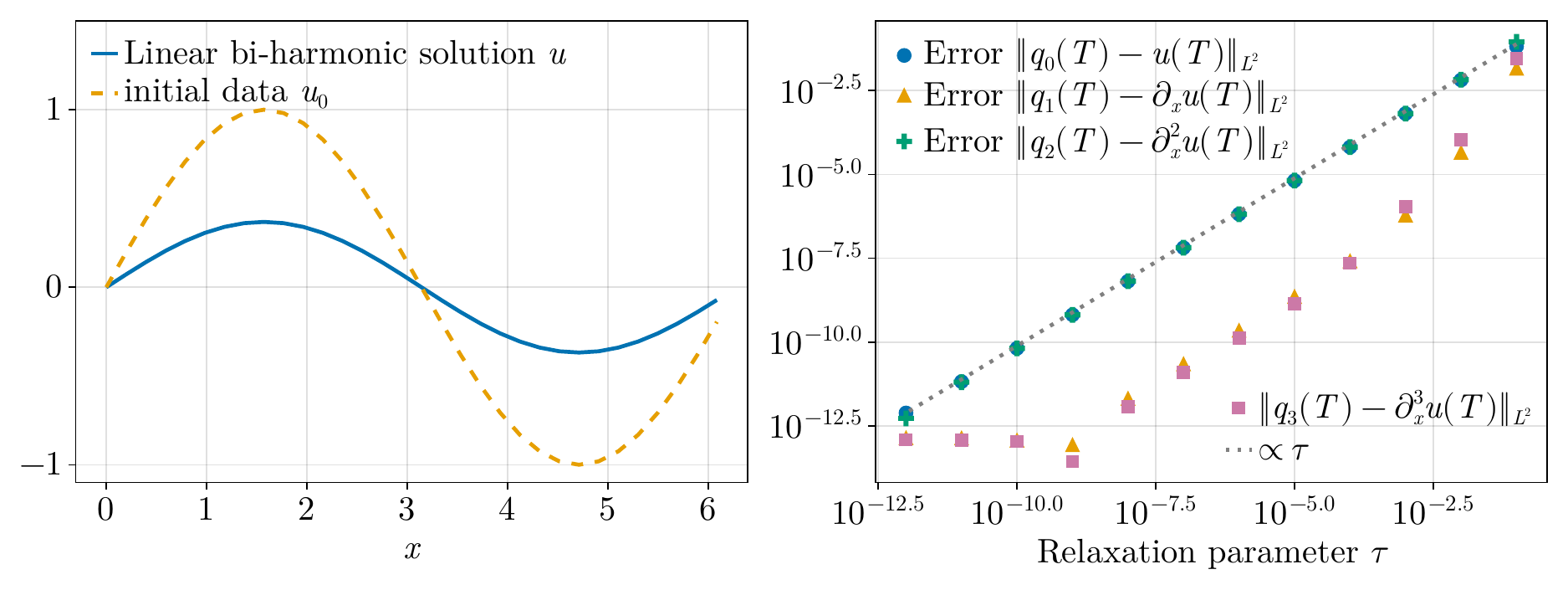}
  \caption{Convergence of the hyperbolic approximation \eqref{eq:spatial_even_hyperbolic}
           to the linear bi-harmonic equation \eqref{eq:linear_biharmonic}.
           The left plot shows the numerical solution at the final time
           $T = 1$ and the initial condition.
           The right plot shows the convergence of the discrete $L^2$ error at the final time $T$
           as a function of the hyperbolic relaxation parameter $\tau$.}
  \label{fig:linear_biharmonic_convergence}
\end{figure}

The results obtained with $2^{5}$ grid points, third-order accurate
upwind operators, and time step size $\Delta t = 0.01$ are shown in
Figure~\ref{fig:linear_biharmonic_convergence}. We observe that the
hyperbolic approximation converges to the solution of the linear
bi-harmonic equation with the expected order $O(\tau)$ predicted
by Theorem~\ref{thm:convergence_spatial_even}.
The even-order derivative-approximation converges with the same order
and the odd-order derivative approximations appear to converge even faster
until $\tau$ is so small that conditioning and rounding errors become an issue.

\subsection{Kuramoto-Sivashinsky equation}

Consider a (sufficiently regular) solution $u$ to the Kuramoto-Sivashinsky equation
\cite{kuramoto1980instability,sivashinsky1977nonlinear,tadmor1986wellposedness}
\begin{equation}
\label{eq:KS_limit}
  \partial_t u + \partial_x \frac{u^2}{2} + \partial_x^2 u + \partial_x^4 u = 0.
\end{equation}
Ketcheson and Biswas \cite{ketcheson2025approximation} proposed the hyperbolic approximation
\begin{equation}
\begin{aligned}
  \partial_t q_0 + \partial_x (q_0^2 / 2) + \partial_x q_1 + \partial_x q_3 &= 0,\\
  \tau \partial_t q_1 - (\partial_x q_2 - q_3) &= 0,\\
  \tau \partial_t q_2 - (\partial_x q_1 - q_2) &= 0,\\
  \tau \partial_t q_3 + (\partial_x q_0 - q_1) &= 0.
\end{aligned}
\end{equation}
If we set $\eta(q):= \int q_0^2 + \tau q_1^2 + \tau q_2^2 + \tau q_3^2 \dif x$, solutions of this system satisfy
\[
 \frac{\dif \eta(q)}{\dif t}= \int_\Omega \left( - q_0 \partial_x q_1 - q_2^2 \right) \dif x,
\]
and it is not clear how to control the $q_0 \partial_x q_1$ term on the right-hand side.

Thus, we choose the slightly different hyperbolic approximation
\begin{equation}\label{eq:KS_hyperbolic}
\begin{aligned}
  \partial_t q_0 + \partial_x (q_0^2 / 2) + q_2 + \partial_x q_3 &= 0,\\
  \tau \partial_t q_1 - (\partial_x q_2 - q_3) &= 0,\\
  \tau \partial_t q_2 - (\partial_x q_1 - q_2) &= 0,\\
  \tau \partial_t q_3 + (\partial_x q_0 - q_1) &= 0,
\end{aligned}
\end{equation}
where we use $q_2$ instead of $\partial_x q_1$ to approximate $\partial_x^2 u$
in the first equation. Note that this is different from our approach to handle
the additional lower-order dispersive term in the Kawahara equation
(see Section~\ref{sec:kawahara}).
For solutions of \eqref{eq:KS_hyperbolic},
\[
 \frac{\dif \eta(q)}{\dif t}= \int_\Omega - q_0 q_2 - q_2^2 \dif x
 \leq \frac14\int_\Omega  q_0^2 \dif x,
\]
which means that the energy might grow but in a controlled --- at most exponential --- way, which is sufficient for our relative energy framework.

The system \eqref{eq:KS_hyperbolic} is indeed hyperbolic with flux Jacobian
\begin{equation}
  \begin{pmatrix}
    q_0 & 0 & 0 & 1 \\
    0 & 0 & -\tau^{-1} & 0 \\
    0 & -\tau^{-1} & 0 & 0 \\
    \tau^{-1} & 0 & 0 & 0
  \end{pmatrix},
\end{equation}
which is diagonalizable with the eigenvalues and eigenvectors given by
\begin{equation}
\begin{gathered}
  \lambda_1 = \frac{-1}{\tau}, \quad
  \lambda_2 = \frac{1}{\tau}, \quad
  \lambda_3 = \frac{q_0 - \sqrt{q_0^2 + 4 / \tau}}{2}, \quad
  \lambda_4 = \frac{q_0 + \sqrt{q_0^2 + 4 / \tau}}{2}, \\
  v_1 = \begin{pmatrix} 0 \\ 1 \\ 1 \\ 0 \end{pmatrix}, \quad
  v_2 = \begin{pmatrix} 0 \\ -1 \\ 1 \\ 0 \end{pmatrix}, \quad
  v_3 = \begin{pmatrix} \tau \lambda_3 \\ 0 \\ 0 \\ 1 \end{pmatrix}, \quad
  v_4 = \begin{pmatrix} \tau \lambda_4 \\ 0 \\ 0 \\ 1 \end{pmatrix}.
\end{gathered}
\end{equation}

Given a sufficiently smooth solution $u$ to \eqref{eq:KS_limit},
we construct an approximate solution $\bar q$ to \eqref{eq:KS_hyperbolic} that satisfies
\begin{equation}
\label{eq:KS_res}
\begin{aligned}
  \partial_t \bar q_0 + \partial_x (\bar q_0^2 / 2) + \bar q_2 + \partial_x \bar q_3 &= \tau R,\\
  \tau \partial_t \bar q_1 - (\partial_x \bar q_2 - \bar q_3) &= 0,\\
  \tau \partial_t \bar q_2 - (\partial_x \bar q_1 - \bar q_2) &= 0,\\
  \tau \partial_t \bar q_3 + (\partial_x \bar q_0 - \bar q_1) &= 0,
\end{aligned}
\end{equation}
where the residual $R$ depends on $u$ and its derivatives up to order $4$
and is uniformly bounded for $\tau \to 0$.
We can do this by defining $\bar q_2 = \partial_{xx} u$. Then we define $\bar q_1$ such that \eqref{eq:KS_res}$_3$ is satisfied, i.e. $\bar q_1=\partial_x u + \tau \partial_{tx} u$. Then we define $\bar q_3$ such that \eqref{eq:KS_res}$_2$ is satisfied, i.e., $\bar q_3 = \partial_{xxx} u - \tau \partial_{tx} u - \tau^2 \partial_{ttx} u$.  Finally, we define $\bar q_0$ such that \eqref{eq:KS_res}$_4$ is satisfied, $\bar q_0 = u + \tau \partial_{t} u - \tau (\partial_{txx} u - \tau \partial_{tt} u - \tau^2 \partial_{ttt} u)$.

We define the relative energy
\begin{equation}
  \eta(q, \bar q)
  =
  \int \left(
    \frac12 (q_0 - \bar q_0)^2 + \frac{\tau}{2} \sum_{j=1}^3 (q_j - \bar q_j)^2
  \right) \dif x.
\end{equation}
Its time evolution can be computed setting $f(u) = u^2 / 2$ as
\begin{align*}
  \frac{\dif}{\dif t} \eta(q,\bar q)
  &=
  \int \left( (q_0- \bar q_0) \partial_t (q_0- \bar q_0) + \sum_{j=1}^{3} \tau (q_j- \bar q_j) \partial_t (q_j- \bar q_j) \right) \dif x
  \\
  &= -\int \biggl(
    (q_0- \bar q_0) \partial_x \big(f(q_0) - f(\bar q_0)\big) +(q_0- \bar q_0)\tau R - (q_0- \bar q_0) (q_2- \bar q_2) - (q_0- \bar q_0) \partial_x (q_3- \bar q_3) \\
    &\qquad + (q_1- \bar q_1)\partial_x (q_2- \bar q_2) - (q_1- \bar q_1) (q_3- \bar q_3) + (q_2- \bar q_2) \partial_x(q_1- \bar q_1) -(q_2- \bar q_2)^2 \\
    & \qquad - (q_3- \bar q_3) \partial_x (q_0 - \bar q_0)   + (q_1- \bar q_1) (q_3- \bar q_3)\biggr)\dif x\\
    &= -\int \biggl(
    (q_0- \bar q_0) \partial_x \big(f(q_0) - f(\bar q_0)\big) +(q_0- \bar q_0)\tau R - (q_0- \bar q_0) (q_2- \bar q_2) - (q_2- \bar q_2)^2 \biggr) \dif x\\
    &\leq \|\partial_x \bar q_0\|_{\infty}  \int \bigl|f(q_0)- f(\bar q_0)- f'(\bar q_0)(q_0- \bar q_0) \bigr|\dif x  +\int \biggl( (q_0- \bar q_0)^2 + \frac{\tau^2}{2} R^2\biggr) \dif x\\
    & \leq ( C \|\partial_x \bar q_0\|_{\infty} + 2)\eta(q, \bar q) + \tau^2 \| R\|_{L^2}^2.
\end{align*}
Thus, we can conclude
\begin{equation}
  \| u - q_0 \|_{L^\infty(0,T, L^2(\Omega))}
  + \sum_{j=1}^{3} \sqrt{\tau} \| \partial_x^j u - q_j \|_{L^\infty(0,T, L^2(\Omega))}
  = \mathcal{O}(\tau)
\end{equation}
by Gronwall's lemma. This proves
\begin{theorem}
\label{thm:convergence_kuramoto_sivashinsky}
   Let $T>0$ and let $u \in W^{4,\infty}((0,T) \times \Omega)$ be a solution to
  \eqref{eq:spatial_odd_limit} with initial data $u_0 \in W^{4,\infty}(\Omega)$.
  Let, for each $\tau>0$, $q$ be an entropy solution to \eqref{eq:spatial_odd_hyperbolic}
  with $q|_{t=0}= (u_0, \partial_x u_0 ,\dots, \partial_x^{m-1} u_0 )$. Then
  \begin{equation}
    \| u - q_0 \|_{L^\infty(0,T, L^2(\Omega))}
    + \sum_{j=1}^{3} \sqrt{\tau} \| \partial_x^j u - q_j \|_{L^\infty(0,T, L^2(\Omega))}
    =
    \mathcal{O}(\tau).
  \end{equation}
\end{theorem}

\begin{remark}
 Note that an analogous statement to Remark \ref{rem:convergence_kawahara} holds in this case.
\end{remark}

\subsubsection{Numerical demonstration}

We discretize the Kuramoto-Sivashinsky equation \eqref{eq:KS_limit} as
\begin{equation}
  \partial_t \vec{u}
  + \frac{1}{3} \left(
    \vec{u} D_1 \vec{u} + D_1 \vec{u}^2
  \right)
  + D_+ D_- \vec{u} + D_+ D_- D_+ D_- \vec{u}
  =
  \vec{0}
\end{equation}
and the hyperbolic approximation \eqref{eq:KS_hyperbolic} as
\begin{equation}
  \partial_t
  \begin{pmatrix}
    \vec{q_0} \\
    \vec{q_1} \\
    \vec{q_2} \\
    \vec{q_3}
  \end{pmatrix}
  + \begin{pmatrix}
      \frac{1}{3} (\vec{q_0} D_1 \vec{q_0} + D_1 \vec{q_0}^{\!\!2}) \\
      \vec{0} \\
      \vec{0} \\
      \vec{0}
  \end{pmatrix}
  - \begin{pmatrix}
    \vec{q_2} + D_+ \vec{q_3} \\
    \tau^{-1} (D_- \vec{q_2} - \vec{q_3}) \\
    \tau^{-1} (D_+ \vec{q_1} - \vec{q_2}) \\
    \tau^{-1} (-D_- \vec{q_0} + \vec{q_1})
  \end{pmatrix}
  =
  \begin{pmatrix}
    \vec{0} \\
    \vec{0} \\
    \vec{0} \\
    \vec{0}
  \end{pmatrix}.
\end{equation}
We treat the nonlinear part explicitly and linear part implicitly.

\begin{figure}[htb]
  \includegraphics[width=\textwidth]{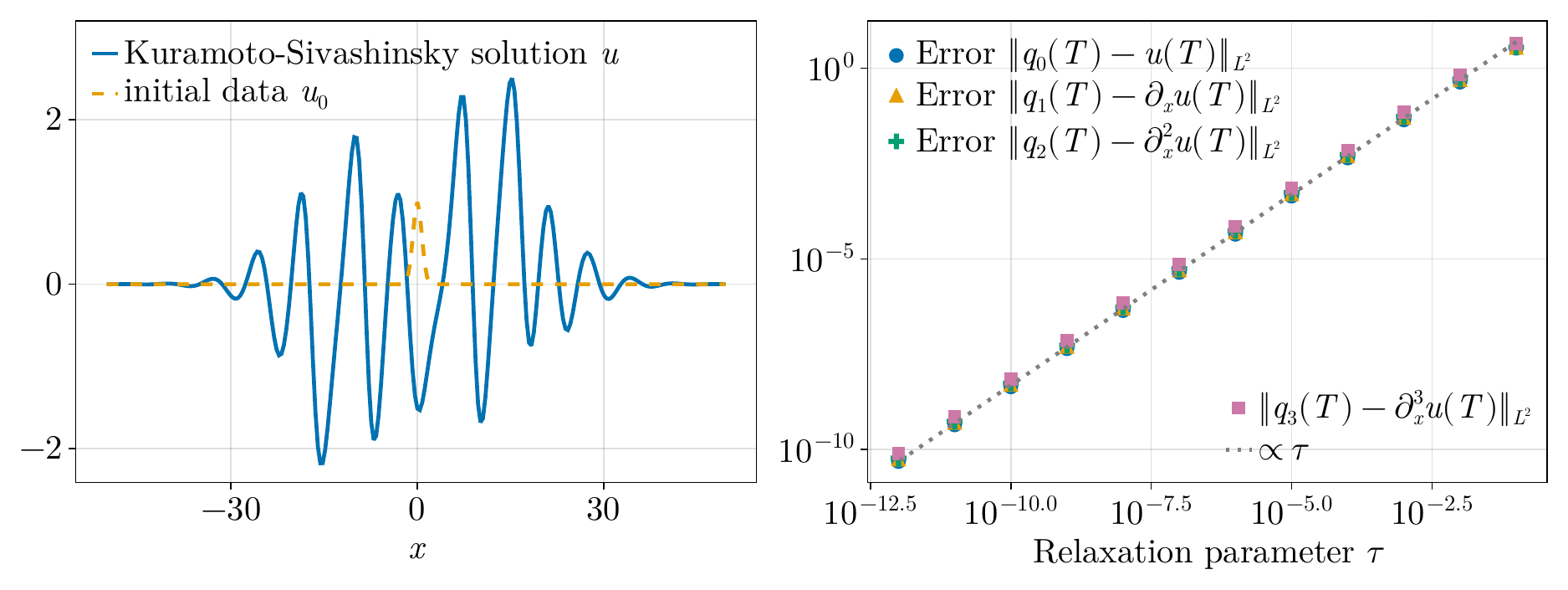}
  \caption{Convergence of the hyperbolic approximation
           \eqref{eq:KS_hyperbolic}
           to the Kuramoto-Sivashinsky equation \eqref{eq:KS_limit}.
           The left plot shows the numerical solution at the final time
           $T = 20$ and the initial condition.
           The right plot shows the convergence of the discrete $L^2$ error at the final time $T$
           as a function of the hyperbolic relaxation parameter $\tau$.}
  \label{fig:kuramoto_sivashinsky_convergence}
\end{figure}

To demonstrate the convergence of the hyperbolic approximation
\eqref{eq:KS_hyperbolic} to the Kuramoto-Sivashinsky equation
\eqref{eq:KS_limit}, we use the initial condition
\begin{equation}
  u(0, x) = \exp\bigl(-x^2\bigr)
\end{equation}
in the domain $\Omega = [-50, 50]$ with periodic boundary conditions.
We use $2^{8}$ grid points, seventh-order accurate upwind operators,
and time step size $\Delta t = 0.1$.
The numerical solution of the Kuramoto-Sivashinsky equation as well as
the discrete $L^2$ error are shown in Figure~\ref{fig:kuramoto_sivashinsky_convergence}.
As expected, the hyperbolic approximation converges to the Kuramoto-Sivashinsky
solution with the expected order $O(\tau)$ predicted by
Theorem~\ref{thm:convergence_kuramoto_sivashinsky}.
Again, the derivative approximations converge with the same order.

\section{Summary and discussion}
\label{sec:summary}

We have proven the convergence of hyperbolic approximations to higher-order PDEs for several classes of equations.
Since our analysis is based on the relative energy/entropy method, it requires smooth solutions of the original PDEs but only weak (entropy) solutions of the hyperbolizations.
Interesting extensions that will require new techniques include convergence results for non-smooth solutions.

Our analysis is supported by numerical results using structure-preserving discretizations based on summation-by-parts operators in space and additive Runge-Kutta methods in time.
The numerical results suggest that the derivative approximations $q_i$, $i \ge 1$, converge with the same order as the leading-order approximation $q_0$ approximating the solution $u$ of the original PDE.
This is an even stronger result that is not fully supported by our theoretical analysis.
Thus, we conjecture that the analysis can be refined to show the improved convergence rates observed in our numerical experiments (and the ones in \cite{biswas2025traveling,bleecke2025asymptotic}).

To keep the focus of our manuscript clear, we have not discussed the details of the numerical methods.
However, we would like to stress that all of them preserve the required structures discretely.
While the numerical results suggest that the methods are indeed asymptotic-preserving, mimicking the convergence rates of the continuous theory, we have not proven this rigorously in this manuscript to keep the focus on the convergence of the hyperbolic approximations.
All of these structure-preserving properties of the discretizations can be analyzed similar to \cite{biswas2025traveling,bleecke2025asymptotic}.

\appendix

\section*{Acknowledgments}

The authors acknowledge funding by the Deutsche Forschungsgemeinschaft (DFG, German Research Foundation) within \emph{SPP 2410 Hyperbolic Balance Laws in Fluid Mechanics: Complexity, Scales, Randomness (CoScaRa)}, project numbers  525877563(JG) and 526031774 (HR).
JG also acknowledges support by the DFG via grant TRR 154 (Mathematical modelling, simulation and optimization using the example of gas networks), subproject C05 (Project 239904186).
HR additionally acknowledges the DFG individual research grants 513301895 and 528753982.

\printbibliography

\end{document}